\documentclass[11pt]{article}
\usepackage[T1]{fontenc}
\usepackage[utf8]{inputenc}
\usepackage{lmodern,amsmath,amssymb,amsthm,mathtools,booktabs,array,graphicx}
\usepackage[margin=1in]{geometry}
\usepackage{microtype}
\usepackage[colorlinks=true,linkcolor=blue,citecolor=blue,urlcolor=blue]{hyperref}
\usepackage{enumitem}
\setlist{nosep}
\newtheorem{theorem}{Theorem}[section]
\newtheorem{proposition}[theorem]{Proposition}
\newtheorem{lemma}[theorem]{Lemma}

\theoremstyle{remark}
\newcommand{\E}{\mathbb E}
\newcommand{\R}{\mathbb R}
\newcommand{\C}{\mathbb C}
\newcommand{\Var}{\operatorname{Var}}
\newcommand{\Cov}{\operatorname{Cov}}
\newcommand{\dd}{\,\mathrm d}

\newcommand{\clip}{\operatorname{clip}}
\newcommand{\Lip}{\operatorname{Lip}}
\newcommand{\norm}[1]{\left\lVert #1\right\rVert}
\DeclareMathOperator{\Rea}{Re}
\DeclareMathOperator{\Ima}{Im}
\title{RH-Conditional Positivity in a Theta-Kernel Langevin Model}
\author{U Jin Choi\\KAIST}
\date{November 2026}
\begin{document}
\maketitle
\begin{abstract}
A classical Brownian-bridge representation of Riemann's xi function defines an even, rapidly confining probability density and a reversible Langevin model. Our principal physical observable is the transient variance after a quadratic-potential quench; we also analyze a tilted bridge-range response and a sixth-order cumulant gap. For a specified passive harmonic trap with bounded delayed residual feedback, we derive a non-vacuous finite-time variance-increase bound under explicit calibration and sampling assumptions. An elementary global curvature estimate and de Bruijn's strip-contraction theorem give a direct proof of a strict susceptibility bound for every nonnegative quadratic tilt, without assuming the Riemann hypothesis. The inequality itself is a classical GHS consequence, valid more generally for every real quadratic tilt. Reversibility then gives explicit finite-time response and deterministic-control action bounds (not thermodynamic energy). The same curvature estimate gives monotone variance relaxation after an outward quadratic quench, with an exponential deficit bound. At the undeformed density, all alternating zero-field cumulant signs and the displayed low-order cumulant gap are also unconditional; neither is an RH test. The genuinely stronger RH-conditional structure is a positive Stieltjes representation and positive semidefiniteness of divided-difference matrices of every size. At the undeformed density that entire matrix hierarchy is equivalent to RH. We separate these exact statements from finite-sample confidence bounds, numerical discretization evidence, and theoretical physical implementations. The work supplies neither a proof of RH nor a laboratory realization.
\end{abstract}

\section{Introduction}
The connection between Brownian functionals and the Riemann xi function is classical. In particular, the scaled range of a unit Brownian bridge has Mellin transform $2\xi$; Biane, Pitman, and Yor review this identity and its probabilistic setting \cite{BPY}. It yields a positive probability model before any hypothesis about zeta zeros is imposed. Exponential tilting of its logarithm turns derivatives of $\log\xi$ into ordinary means, variances, and cumulants.

Our purpose is to distinguish consequences that can already be proved from consequences that genuinely require the full real-zero property. This distinction matters for experimental interpretation. A scalar response inequality, a positive low-order cumulant determinant, or a finite collection of successful numerical tests cannot establish a global statement about all complex zeros. Conversely, a useful physical inequality need not be described as conditional on RH when it follows from weaker, established information.

The mathematical ingredients have substantial precedents. De Bruijn's deformation and Newman's work relate Gaussian kernel tilts to real zeros \cite{deBruijn,Newman1976,NewmanWu}; Rodgers and Tao prove the nonnegativity of the de Bruijn--Newman constant \cite{RodgersTao}. Newman derived correlation inequalities from Lee--Yang zeros \cite{Newman1975}. Reversible response formulas belong to the fluctuation--dissipation framework \cite{Kubo}, while matrix positivity is governed by L\"owner's theorem \cite{Sparr}. We do not claim priority for these principles, the bridge identity, the heat deformation, or the general connection between partition-function zeros and cumulants. Our contribution is the quantitative analysis of this specified reversible diffusion and its finite-resource observation schemes. In particular, the susceptibility inequality is not claimed as new: Newman's theta-kernel result \cite{Newman1991}, combined with the GHS theorem of Ellis, Monroe, and Newman \cite{EMN}, already gives a stronger all-real-parameter conclusion. The 2026 preprint of Min and Kang \cite{MinKang} explicitly formulates the RH--Stieltjes and Hankel equivalence; it is cited as a preprint, not as peer-reviewed validation. Suzuki \cite{Suzuki} studies a closely related square-variable Stieltjes representation, and Pascoe \cite{Pascoe} gives an RH characterization through matrix monotonicity of a logarithmic derivative, in a different real-variable normalization. The self-contained matrix argument below is an application of classical L\"owner theory, not a claimed new RH reformulation.

High-order cumulants have previously been used to characterize Lee--Yang zeros and dynamical phase transitions by Flindt and Garrahan \cite{FlindtGarrahan}. For a broader survey of physical approaches to RH, see Schumayer and Hutchinson \cite{SchumayerHutchinson}. Physical access to complex partition-function zeros also has direct precedents: Peng and collaborators reported an experimental observation of Lee--Yang zeros \cite{Peng}, and Wei and collaborators studied the Riemann-hypothesis connection in dynamical quantum phase transitions \cite{Wei}. Those works concern different systems and do not constitute a realization of the overdamped theta-kernel potential proposed here. We claim neither the first physical interpretation of xi nor the first experimental route to complex zeros.

\paragraph{The three observation routes.}
Route C is the principal physical analysis: it perturbs the confinement quadratically and measures the transient variance. Route A measures a nonlinear mean through a weighted Brownian-bridge range. Route B estimates a zero-field cumulant gap, with explicit treatment of its sampling covariance. All three are well-defined without RH. They interrogate different features of the same probability family and have different error budgets.

\section{Methods: The probability family and its normalizations}
\paragraph{AI use in the research method.} Under the author's direction, AI (ChatGPT, OpenAI; version not recorded) assisted with mathematical exploration and preparation of numerical scripts. Checks performed during preparation included symbolic-algebra checks, rerunning the principal numerical calculations, and further AI-assisted scrutiny of proofs; independent author verification is still required before submission. The author retains responsibility for accuracy and the final content.
Let $B$ be a standard Brownian bridge on $[0,1]$, with covariance
$\E B_sB_t=\min(s,t)-st$, and put
\begin{equation}\label{eq:bridge}
 R_B=\sup_{0\leq t\leq1}B_t-\inf_{0\leq t\leq1}B_t,
 \qquad Y=\sqrt{2/\pi}\,R_B,\qquad X=\log Y.
\end{equation}
Write $P$ for this bridge law. The classical identity is
\begin{equation}\label{eq:mellin}
 \E_P Y^s=2\xi(s),\qquad
 \xi(s)=\tfrac12s(s-1)\pi^{-s/2}\Gamma(s/2)\zeta(s),
 \qquad s\in\C,
\end{equation}
with removable singularities filled in \cite{BPY}. All complex moments exist. Define
\begin{equation}\label{eq:sizebias}
 \frac{\dd Q_0}{\dd P}=\frac{Y^{1/2}}{2\xi(1/2)}.
\end{equation}
Under $Q_0$, $X$ has an even, positive, smooth density $q_0$ and
\begin{equation}\label{eq:mgf0}
 M_0(h)=\E_{Q_0}e^{hX}
       =\frac{\xi(1/2+h)}{\xi(1/2)}.
\end{equation}
Evenness follows directly from $\xi(s)=\xi(1-s)$ and uniqueness of the moment generating function; it is not an RH assumption.

For a convenient kernel normalization let, for $x\geq0$,
\begin{equation}\label{eq:kernel}
 K(x)=\sum_{n\geq1}e^{x/2}R_n(2R_n-3)e^{-R_n},
 \qquad R_n=\pi n^2e^{2x},\qquad K(-x)=K(x).
\end{equation}
Then $q_0=K/\int_\R K$, and $\int_\R K=\xi(1/2)/2$.
The theta functional equation makes the even extension smooth at zero.
All terms in \eqref{eq:kernel} are positive on $x\geq0$. The tails are double exponential, so that every polynomial, linear exponential, and fixed quadratic exponential is integrable.

To make the extension at zero explicit, let
$\vartheta(t)=\sum_{n\in\mathbb Z}e^{-\pi n^2t}$ for $t>0$.
With $t=e^{2x}$ the displayed kernel series, evaluated for any real
$x$, equals
\[
 \mathcal K(x)=t^{1/4}
 \left(t^2\vartheta''(t)+\frac32t\vartheta'(t)\right).
\]
The series and all its derivatives converge uniformly on compact
$x$-intervals, by a bound of the form $C n^d e^{-c n^2}$.
The classical theta functional equation is
$\vartheta(t)=t^{-1/2}\vartheta(1/t)$.
Differentiating it twice, with $s=1/t$, gives
\[
 t^{1/4}\left(t^2\vartheta''(t)+\frac32t\vartheta'(t)\right)
 =s^{1/4}\left(s^2\vartheta''(s)+\frac32s\vartheta'(s)\right).
\]
For example the differentiated formulas are
$\vartheta'(t)=-\tfrac12t^{-3/2}\vartheta(s)
-t^{-5/2}\vartheta'(s)$ and
$\vartheta''(t)=\tfrac34t^{-5/2}\vartheta(s)
+3t^{-7/2}\vartheta'(s)+t^{-9/2}\vartheta''(s)$;
the undifferentiated terms cancel in the indicated combination.
Thus $\mathcal K(-x)=\mathcal K(x)$.  The even extension in
\eqref{eq:kernel} is this smooth function, rather than a merely
piecewise-smooth definition.

For $a,h\in\R$ define
\begin{align}
 Z(a,h)&=\int_\R e^{ax^2+hx}q_0(x)\dd x,
 &q_{a,h}(x)&=Z(a,h)^{-1}e^{ax^2+hx}q_0(x),\\
 M_a(h)&=Z(a,h)/Z(a,0),
 &m_a(h)&=\partial_h\log M_a(h),\\
 \chi_a(h)&=\partial_hm_a(h)=\Var_{q_{a,h}}X,
 &v_a&=\chi_a(0).
\end{align}
We abbreviate $q_a=q_{a,0}$ and
$\phi_a(z)=\int e^{izx}q_a(x)\dd x$.
Both transforms are entire and $M_a(h)=\phi_a(-ih)=\phi_a(ih)$.
To avoid a common factor-of-four error, in the convention
\begin{equation}\label{eq:Hnormal}
 H_t(z)=\int_0^\infty e^{tu^2}\Phi(u)\cos(zu)\dd u,
 \qquad \Phi(u)=K(2u),
\end{equation}
one has exactly
\begin{equation}\label{eq:scaling}
 \phi_a(z)=\frac{H_{4a}(2z)}{H_{4a}(0)},\qquad
 H_0(z)=\frac18\xi(1/2+iz/2).
\end{equation}
This is the convention used in \cite{RodgersTao}.

\section{An unconditional response theorem}
\begin{theorem}\label{thm:static}
For every $a\geq0$ and $h\neq0$,
\begin{equation}\label{eq:chain}
 0<\chi_a(h)<\frac{m_a(h)}h<v_a.
\end{equation}
At $h=0$ the last two comparisons become equalities by continuous extension. No Riemann hypothesis or finite-height verification of RH is needed for this theorem.
\end{theorem}

\paragraph{Classical GHS scope.}
Two external theorems give a stronger parameter range than the
independent proof below.  Newman \cite{Newman1991} proves that
$V'$, for the even theta-kernel potential $V=-\log q_0$, is convex
on $[0,\infty)$; a change of spatial scale and an additive
normalization constant do not affect this property.
For each real $a$ the potential $V_a(x)=V(x)-ax^2$ is still even,
its exponential and all its linear tilts are integrable, and
$V_a'=V'-2ax$ is convex on the positive half-line, because
subtracting an affine function preserves convexity.
The one-site case of the GHS theorem of Ellis, Monroe, and Newman
\cite{EMN} therefore gives
$\partial_h^3\log M_a(h)=\chi_a'(h)\leq0$ for $h\geq0$.
These are cited classical results, not conclusions derived from
RH or proved by the elementary curvature estimate.

For $h>0$, integrating $m_a'=\chi_a$ gives
$m_a(h)/h=h^{-1}\int_0^h\chi_a(s)\dd s$.
Monotonicity implies
$\chi_a(h)\leq m_a(h)/h\leq\chi_a(0)$.
If equality held in either comparison at a positive $h$, continuity
and monotonicity would force $\chi_a$ to be constant throughout
$[0,h]$.  The moment-generating function is entire and positive on
the real axis, so $\chi_a=(\log M_a)''$ is real analytic there.
The real-analytic identity theorem makes it constant on all of
$\R$.  Since $M_a(0)=1$ and $m_a(0)=0$, integration would give
$M_a(h)=\exp(v_ah^2/2)$, the moment-generating function of a
centered Gaussian.  Uniqueness of the moment-generating function
would make $q_a$ Gaussian, contradicting its double-exponential
tails.  Hence both comparisons are strict.  Evenness extends them
to negative $h$, and nondegeneracy supplies $\chi_a(h)>0$.
Thus the strict chain holds for every real $a$ by these classical
theorems.  The direct proof that follows deliberately treats
$a\geq0$ and does not claim a new susceptibility inequality.

We give that direct proof in three steps. The first also supplies the contraction constant used later.
\begin{lemma}\label{lem:curvature}
For every $x\in\R$,
\begin{equation}\label{eq:curvature}
 - (\log q_0)''(x)>18,
 \qquad b_0'(x)<-9,
 \qquad b_0=\tfrac12(\log q_0)'.
\end{equation}
\end{lemma}
\begin{proof}
For $x\geq0$ put $R=\pi e^{2x}$ and
$f_n(x)=e^{x/2}n^2R(2n^2R-3)e^{-n^2R}$.  These functions are
positive.  On every compact $x$-interval, each differentiated summand
is bounded by a polynomial in $n$ times $e^{-c n^2}$ for a positive
constant $c$ depending on the interval.  Thus all the differentiations
of the series below are locally uniformly convergent.  With
$w_n=f_n/K$, direct differentiation gives
\begin{align*}
 \beta_n=(\log f_n)'&=\frac92+\frac6{2n^2R-3}-2n^2R,\\
 k_n=-(\log f_n)''&=4n^2R+\frac{24n^2R}{(2n^2R-3)^2}.
\end{align*}
In particular, $w_n'=w_n(\beta_n-\sum_jw_j\beta_j)$, so
\[
 (\log K)''=-\sum_nw_nk_n+
 \sum_nw_n\beta_n^2-\left(\sum_nw_n\beta_n\right)^2.
\]
We bound the curvature and variance terms separately.

Write $k(r)=4r+24r/(2r-3)^2$.  Since $n^2R\geq\pi>3$,
it suffices to consider $r\geq3$.  For $y=r-3\geq0$ the exact
algebraic identity is
\begin{align*}
 (k(r)-19)(2r-3)^2
 &=16r^3-124r^2+288r-171\\
 &=9-24y+20y^2+16y^3\\
 &=20(y-3/5)^2+9/5+16y^3>0.
\end{align*}
Thus $k_n>19$ for every $n$.

For any constant $b$, the variance-minimization identity gives
$\Var_w(\beta_n)\leq\sum_nw_n(\beta_n-b)^2$.
Use $b=\beta_1$, $K\geq f_1$, $d=n^2-1$, and
$C=6/(2\pi-3)<2$.  Then
\[
 0\leq\beta_1-\beta_n
 =2dR+\frac6{2R-3}-\frac6{2n^2R-3}\leq2dR+C,
\]
and therefore
\[
 \Var_w(\beta_n)
 \leq\sum_{n\geq2}A_n(R)(2dR+C)^2e^{-dR},
 \qquad
 A_n(R)=\frac{n^2(2n^2R-3)}{2R-3}.
\]
For $n\geq2$ and $R\geq\pi$,
\[
 A_n'(R)=-\frac{6n^2(n^2-1)}{(2R-3)^2}<0,
 \qquad
 \frac{\dd}{\dd R}\log\{(2dR+C)^2e^{-dR}\}
 =d\left(\frac4{2dR+C}-1\right)<0.
\]
Each summand is consequently at most its value at $R=\pi$.
The elementary bounds $157/50<\pi<22/7$ imply
$A_2(\pi)=16+36/(2\pi-3)<27$ and $6\pi+C<21$.
Also $e^3>20$ follows from its Taylor series, and hence
\[
 e^{3\pi}>e^{9+21/50}
 >20^3\left(1+\frac{21}{50}
              +\frac{(21/50)^2}{2}
              +\frac{(21/50)^3}{6}\right)>12000.
\]
The $n=2$ contribution is less than $27\cdot21^2/12000=0.99225$.
For $n\geq3$ we have
\[
 A_n(\pi)=n^4+\frac{3n^2(n^2-1)}{2\pi-3}<2n^4,
 \qquad 2(n^2-1)\pi+C<7n^2.
\]
Thus the remaining contribution is at most
$\sum_{n\geq3}98n^8e^{-(n^2-1)\pi}$.
The first term is less than $98\cdot3^8/20^8$ and the ratio of
successive terms is at most
\[
 (4/3)^8e^{-7\pi}<\frac{10}{20^7}.
\]
Summing this geometric majorant gives
\[
 \sum_{n\geq3}98n^8e^{-(n^2-1)\pi}
 <\frac{98\cdot3^8/20^8}{1-10/20^7}<0.000025117.
\]
Consequently $\Var_w(\beta_n)<1$ and
$-(\log K)''>19-1=18$ on $[0,\infty)$.
The theta functional equation gives the smooth even extension of $K$;
thus the same estimate holds on the negative half-line and at zero.
Since normalization does not change logarithmic derivatives,
$-(\log q_0)''>18$ and $b_0'=(\log q_0)''/2<-9$.
\end{proof}

\begin{lemma}\label{lem:zeros}
For $0\leq a\leq1/8$, every zero $z$ of $\phi_a$ satisfies
\begin{equation}\label{eq:rectangle}
 |\Ima z|\leq\sqrt{1/4-2a},\qquad |\Rea z|>2.
\end{equation}
For $a\geq1/8$, every zero is real and nonzero.
\end{lemma}
\begin{proof}
We use two classical results, with their normalizations specified.
First, de Bruijn's strip-contraction theorem says that if the Fourier
transform of an even real kernel $k$ has all its zeros in
$|\Ima z|\leq\Delta$, then multiplication of the kernel by
$e^{\lambda x^2/2}$, for $\lambda\geq0$, confines the new zeros to
\[
 |\Ima z|\leq\sqrt{\max(\Delta^2-\lambda,0)}.
\]
The theorem applies to kernels bounded by
$C\exp(-|x|^{2+\epsilon})$ for some $\epsilon>0$;
see \cite[Theorem 7]{NewmanWu}.  Our positive even kernel has
stronger, double-exponential decay.  The critical strip and
\eqref{eq:mgf0} give $\Delta=1/2$ at $a=0$, and the multiplier
in the present notation has $\lambda=2a$.  Normalizing the integral
by its positive value at zero leaves its zeros unchanged.  This proves
the stated strip for $0\leq a\leq1/8$ and real-rootedness for
$a\geq1/8$.  Also $\phi_a(iv)=M_a(-v)>0$ for every real $v$,
so zero and the entire imaginary axis are zero-free.

Second, the one-dimensional Brascamp--Lieb variance inequality states
that a probability density proportional to $e^{-U}$ with $U''\geq\rho>0$
satisfies
$\Var(g)\leq\E[(g')^2/U'']\leq\rho^{-1}\E[(g')^2]$;
see \cite{BrascampLieb}.  For $q_{a,v}$, Lemma~\ref{lem:curvature}
gives $U''>18-2a$.  Applying the inequality to $g(x)=x$ (or first
to smooth cutoffs and then taking the integrable limit) yields
\[
 \chi_a(v)\leq c_a:=\frac1{18-2a}\leq\frac4{71}
 \quad(0\leq a\leq1/8, v\in\R).
\]
Because $m_a(0)=0$ and $m_a'=\chi_a$, integration gives
$|m_a(v)|\leq c_a|v|$.
For $u,v\in\R$, write the complex Fourier transform as a real
exponential tilt:
\[
 \phi_a(u+iv)=M_a(-v)\E_{q_{a,-v}}e^{iuX}.
\]
Taking real parts and using $\cos t\geq1-t^2/2$ gives
\[
 \Rea\phi_a(u+iv)\geq M_a(-v)
 \left[1-\frac{u^2}{2}
  \{\chi_a(-v)+m_a(-v)^2\}\right].
\]
For $|u|\leq2$ and $|v|\leq1/2$ the bracket is at least
\[
 1-2\left(\frac4{71}+\frac4{71^2}\right)>0.
\]
All possible zeros are already in $|v|\leq1/2$, so none has
$|u|\leq2$.  This proves the asserted strict horizontal exclusion.
\end{proof}

\begin{lemma}[Entire growth and logarithmic derivative]\label{lem:analyticproduct}
For each real $a$, there is an entire function $\Psi_a$ of order at most
$1/2$ such that $M_a(h)=\Psi_a(h^2)$ and $\Psi_a(0)=1$.
Index the zeros $z_j$ of $\phi_a$ in $\Rea z>0$, with multiplicity.
Then
\begin{equation}\label{eq:productcomplex}
 \Psi_a(w)=\prod_j\left(1+\frac{w}{z_j^2}\right),\qquad
 M_a(h)=\prod_j\left(1+\frac{h^2}{z_j^2}\right),\qquad
 \sum_j|z_j|^{-2}<\infty.
\end{equation}
There are infinitely many distinct zeros.  Off their poles the series
\begin{equation}\label{eq:meromorphicF}
 \frac{\Psi_a'(w)}{\Psi_a(w)}
 =\sum_j\frac1{w+z_j^2},\qquad
 \left(\frac{\Psi_a'}{\Psi_a}\right)'(w)
 =-\sum_j\frac1{(w+z_j^2)^2}
\end{equation}
converge locally uniformly.
\end{lemma}
\begin{proof}
For $x\geq0$, the kernel series implies
\[
 K(x)\leq D\exp\{(9/2)x-\pi e^{2x}\},\qquad
 D=2\pi^2\sum_{n\geq1}n^4e^{-\pi(n^2-1)}<\infty.
\]
Indeed $R_n(2R_n-3)\leq2R_n^2$ and
$e^{-\pi(n^2-1)e^{2x}}\leq e^{-\pi(n^2-1)}$.
Evenness gives the corresponding bound for negative $x$.
For fixed $a$ choose a finite $C_a$ such that
$ax^2\leq(\pi/4)e^{2x}+C_a$ for $x\geq0$.
This is possible because the exponential dominates the quadratic.
For $|z|\leq r$, the transform and its differentiated integrands
are dominated by integrable double-exponential majorants.  Thus
$\phi_a$ and $M_a$ are entire.  More quantitatively, with
$s=r+9/2$ and a constant $D_a<\infty$,
\begin{align*}
 \max_{|z|\leq r}|\phi_a(z)|
 &\leq D_a\int_0^\infty
       e^{sx-(3\pi/4)e^{2x}}\dd x\\
 &\leq D_a\exp\left\{\sup_{x\geq0}
             \bigl(sx-(\pi/2)e^{2x}\bigr)\right\}
       \int_0^\infty e^{-(\pi/4)e^{2x}}\dd x.
\end{align*}
The supremum, for $s\geq\pi$, equals
$(s/2)\log(s/\pi)-s/2$; for bounded $s$ it is bounded.
It follows that
$\log\max_{|z|\leq r}|\phi_a(z)|=O_a(r\log(r+2))$.
The same bound holds for $M_a(z)=\phi_a(-iz)$.

Since $M_a$ is even, its Taylor series has the form
$M_a(h)=\sum_{n\geq0}c_{2n}h^{2n}$.
The series $\Psi_a(w)=\sum_{n\geq0}c_{2n}w^n$ is entire, and
\[
 \max_{|w|\leq r}|\Psi_a(w)|
 =\max_{|h|\leq\sqrt r}|M_a(h)|,
\]
so $\Psi_a$ has order at most $1/2$.
We now use the classical Hadamard factorization theorem in its
order-less-than-one form: an entire function of order less than one,
nonzero at the origin, is its value at the origin times the genus-zero
product over its zeros, and the sum of the inverse moduli of those
zeros is finite.  There is no nonconstant exponential factor in this
form of the theorem.

There are no imaginary-axis zeros of $\phi_a$, because
$\phi_a(iv)=M_a(-v)>0$.  Its zeros therefore occur in pairs
$\{z,-z\}$, with exactly one member in the right half-plane.
The corresponding zero of $\Psi_a$ is $-z^2$, with the same
multiplicity: the squaring map has nonzero derivative at each
nonzero preimage.  Applying Hadamard's theorem to $\Psi_a$ proves
\eqref{eq:productcomplex}.  If there were only finitely many distinct
zeros, $\Psi_a$, and hence $\phi_a$, would be a polynomial.
But $q_a\in L^1(\R)$, so the Riemann--Lebesgue lemma gives
$\phi_a(t)\to0$ as $|t|\to\infty$ on the real axis.  No nonzero
polynomial has this property, whereas $\phi_a(0)=1$.

Finally, on a compact set $|w|\leq R$ avoiding the poles, all but
finitely many $j$ satisfy $|z_j|^2\geq2R$ and hence
$|w+z_j^2|\geq|z_j|^2/2$.
The series of reciprocals is uniformly dominated by
$2\sum_j|z_j|^{-2}$, and its derivative series by
$4\sum_j|z_j|^{-4}$, both finite after discarding finitely many terms.
Termwise logarithmic differentiation of the normally convergent
product now yields \eqref{eq:meromorphicF}.
\end{proof}

\begin{proof}[Proof of Theorem~\ref{thm:static}]
Put $f_a(y)=m_a(\sqrt y)/\sqrt y$ for $y>0$ and
$f_a(0)=v_a$.  The removable value at zero follows either from
$m_a(0)=0$ and $m_a'(0)=v_a$, or from
$M_a(h)=\Psi_a(h^2)$ in Lemma~\ref{lem:analyticproduct}.
For real $y\geq0$, that lemma gives
\begin{equation}\label{eq:fcomplex}
 f_a(y)=2\sum_j\Rea\frac1{y+z_j^2},\qquad
 f_a'(y)=-2\sum_j\Rea\frac1{(y+z_j^2)^2}.
\end{equation}
The real-part notation does not change the sums: nonreal zeros in
the right half-plane occur together with their conjugates, and all
series converge absolutely.

If $0\leq a\leq1/8$, write $z_j=u+iv$.
Lemma~\ref{lem:zeros} gives $u>2$ and $|v|\leq1/2$.
Set $A=y+u^2-v^2$ and $B=2uv$.  Since $|v|/u<1/4$,
\[
 A\geq u^2-v^2>2u|v|=|B|.
\]
Consequently
\[
 \Rea\frac1{(y+z_j^2)^2}
 =\frac{A^2-B^2}{(A^2+B^2)^2}>0.
\]
If $a\geq1/8$, every $z_j$ is real and the same positivity is
immediate.  There is at least one zero by
Lemma~\ref{lem:analyticproduct}, so $f_a'(y)<0$ for every $y\geq0$.
For $h\neq0$, differentiation of $m_a(h)=h f_a(h^2)$ gives
\[
 \chi_a(h)=f_a(h^2)+2h^2 f_a'(h^2)<f_a(h^2)<f_a(0)=v_a.
\]
The left-hand side is strictly positive because it is the variance
of the positive, nondegenerate density $q_{a,h}$.
This proves every comparison in \eqref{eq:chain}, for either sign
of $h$.  Continuity supplies the equalities at zero.
\end{proof}

\section{The stronger structure supplied by RH}
\begin{theorem}\label{thm:RH}
Assume RH. For each $a\geq0$, write the positive real zeros of $\phi_a$ as $r_j(a)$, with multiplicity. Then
\begin{equation}\label{eq:Stieltjes}
 F_a(y):=\frac{m_a(\sqrt y)}{2\sqrt y}
       =\sum_j\frac1{y+r_j(a)^2},\qquad y\geq0.
\end{equation}
This is a positive Stieltjes transform. For any positive nodes $y_1,\ldots,y_n$, the matrix
\begin{equation}\label{eq:Loewner}
 L_{ij}=\begin{cases}
 \dfrac{F_a(y_i)-F_a(y_j)}{y_j-y_i},&y_i\ne y_j,\\[5pt]
 -F_a'(y_i),&y_i=y_j
 \end{cases}
\end{equation}
is positive semidefinite. It is positive definite for pairwise distinct nodes.
At $a=0$, positivity for every size and every choice of positive nodes is equivalent to RH. The same representation and matrix positivity hold unconditionally for $a\geq0.05$.
\end{theorem}
\begin{proof}
Assume first RH.  Then all zeros of $\phi_0$ are real.
De Bruijn's theorem, applied with initial strip width zero, preserves
real-rootedness under every multiplier $e^{ax^2}$ with $a\geq0$.
There are no zeros at the origin.  Thus the right-half-plane zeros
in Lemma~\ref{lem:analyticproduct} are precisely the positive zeros
$r_j(a)$, counted with multiplicity.  Since
$m_a(h)=2h\Psi_a'(h^2)/\Psi_a(h^2)$, that lemma gives
\eqref{eq:Stieltjes}, including its removable value at $y=0$.

The measure $\mu_a=\sum_j\delta_{r_j(a)^2}$ is locally finite:
an entire function has only finitely many zeros in each compact set,
counted with multiplicity.  Moreover,
\[
 \int_{[0,\infty)}\frac{\dd\mu_a(s)}{1+s}
 =\sum_j\frac1{1+r_j(a)^2}<\infty,
\]
by $\sum_jr_j(a)^{-2}<\infty$.  Therefore the series is a
positive Stieltjes transform in the usual sense.
For distinct nodes, subtracting its two absolutely convergent series
gives
\[
 \frac{F_a(y_i)-F_a(y_j)}{y_j-y_i}
 =\sum_k\frac1{(y_i+r_k^2)(y_j+r_k^2)}.
\]
For equal nodes the same formula follows by the locally uniformly
convergent derivative series in Lemma~\ref{lem:analyticproduct}.
All entries therefore form a Gram matrix:
\[
 c^TLc=\sum_k\left(\sum_{i=1}^n\frac{c_i}{y_i+r_k^2}\right)^2
 \geq0\qquad(c\in\R^n).
\]
The series is finite, since its large-$k$ terms are bounded by a
constant times $r_k^{-4}$.

Suppose now that the nodes are pairwise distinct and $c^TLc=0$.
Every summand is nonnegative, so the rational function
$R(s)=\sum_i c_i/(y_i+s)$ vanishes at every $s=r_k^2$.
There are infinitely many distinct such values by
Lemma~\ref{lem:analyticproduct}.  Multiplication by
$\prod_i(y_i+s)$ gives a polynomial of degree at most $n-1$ with
infinitely many roots, hence the zero polynomial.  Evaluating this
polynomial at $s=-y_i$ gives
$c_i\prod_{j\ne i}(y_j-y_i)=0$.  Distinctness of the nodes implies
$c_i=0$ for every $i$.  This proves positive definiteness without any
assumption of simple zeros.

For the converse fix $a=0$ and assume positivity of all the stated
matrices, of every size.  Set $g(y)=-F_0(y)$ on $(0,\infty)$.
Its ordinary L\"owner matrix has entries
\[
 \frac{g(y_i)-g(y_j)}{y_i-y_j}
 =\frac{F_0(y_i)-F_0(y_j)}{y_j-y_i},\qquad
 g'(y_i)=-F_0'(y_i),
\]
which are exactly the entries in \eqref{eq:Loewner}.
We use the classical L\"owner theorem in the following precise form:
a real differentiable function on an interval whose L\"owner matrices
are positive semidefinite at every size admits a holomorphic Pick
continuation to the upper half-plane, analytic through the interval
and agreeing there with the function \cite{Sparr,Pascoe}.
The hypothesis applies because $F_0$ is real analytic on the positive
axis; indeed $\Psi_0(y)=M_0(\sqrt y)>0$ there.
Let $P$ denote the resulting continuation of $g$.

The quotient $\Psi_0'/\Psi_0$ is a meromorphic function on the whole
plane.  In a neighborhood of each positive real point it equals
$F_0$, so the holomorphic function
\[
 Q(w)=\Psi_0'(w)+P(w)\Psi_0(w),\qquad \Ima w>0,
\]
vanishes on a nonempty open subset of the upper half-plane.
The identity theorem yields $Q\equiv0$ there.  In particular,
$\Psi_0$ cannot have a zero $w_0$ in that half-plane.  To see this
without a removable-singularity assumption, write locally
$\Psi_0(w)=(w-w_0)^m A(w)$ with $m\geq1$ and $A(w_0)\neq0$.
The derivative has a nonzero term of order $m-1$, whereas
$P\Psi_0$ has a zero of order at least $m$, contradicting $Q=0$.
Equivalently, $F_0$ has residue $m$ and $-F_0$ residue $-m$ at
this pole, so multiplicities cannot cause cancellation.

If $z=u+iv$ were a nonreal zero of $\phi_0$ in the right half-plane,
then $u>0$ and $v\neq0$.  Its conjugate is also a zero, and the two
corresponding zeros of $\Psi_0$ are
$-z^2$ and $-\overline z^{\,2}$, whose imaginary parts are
$-2uv$ and $2uv$.  One lies in the upper half-plane, a contradiction.
Imaginary-axis zeros of $\phi_0$ have already been excluded.  Thus
every zero of $\phi_0$ is real, which is RH by \eqref{eq:mgf0}.

Finally, no RH assumption is needed once real-rootedness is known
by another theorem.  Lemma~\ref{lem:zeros} supplies it for
$a\geq1/8$.  The stronger established bound
$\Lambda\leq0.2$ of \cite[Corollary 2]{PlattTrudgian}, together
with the exact parameter relation $t=4a$ in \eqref{eq:scaling},
supplies it for every $a\geq0.05$.  The Stieltjes and Gram arguments
above then apply unchanged.
\end{proof}

\paragraph{What finite measurements can establish.}
A certified negative eigenvalue of one exact matrix would contradict the corresponding real-zero property. A finite list of positive matrices establishes only those necessary inequalities. A negative eigenvalue measured in a physical device challenges RH only if a rigorous uncertainty and model-discrepancy transfer proves negativity for the exact mathematical matrix; otherwise it may indicate device or calibration error. Measurement errors require simultaneous enclosures of all entries; near-coincident nodes amplify error in divided differences. No finite successful experiment verifies the universal quantifiers in Theorem~\ref{thm:RH}. Complete monotonicity by itself is not asserted here to be equivalent to RH.

The connection between real zeros and cumulant inequalities is classical \cite[Proposition~2 and Eqs.~(2.3)--(2.4)]{Newman1975}. Assume RH and let $a\geq0$, so that Theorem~\ref{thm:RH} applies. Define $\kappa_{2j}(a)=\left.\partial_h^{2j}\log M_a(h)\right|_{h=0}$. Here the coefficients follow directly from the product, with all exchanges of series justified as follows.
Let $r_*>0$ be the smallest positive zero and choose
$|h|\leq\eta r_*$ with $0<\eta<1$.
Then
\[
 \sum_k\sum_{j\geq1}\frac{|h|^{2j}}{j r_k^{2j}}
 \leq\frac{|h|^2}{1-\eta^2}\sum_k r_k^{-2}<\infty.
\]
The logarithm normalized by $\log M_a(0)=0$ can consequently be
expanded and rearranged in a neighborhood of zero:
\[
 \log M_a(h)=\sum_{j\geq1}\frac{(-1)^{j-1}}{j}h^{2j}S_j(a),
 \qquad S_j(a)=\sum_k r_k(a)^{-2j}.
\]
Taking the $2j$th derivative at zero gives
\begin{equation}\label{eq:cumulants}
 \kappa_{2j}(a)=2(-1)^{j-1}(2j-1)!S_j(a).
\end{equation}
For $t_k=r_k^{-2}$ and any $c\in\R^n$,
\begin{align*}
 \sum_{i,j=0}^{n-1}c_ic_jS_{i+j+1}
 &=\sum_k t_k\left(\sum_{i=0}^{n-1}c_it_k^i\right)^2\geq0,\\
 \sum_{i,j=0}^{n-1}c_ic_jS_{i+j+2}
 &=\sum_k t_k^2\left(\sum_{i=0}^{n-1}c_it_k^i\right)^2\geq0.
\end{align*}
These sums are finite because $\sum_k t_k<\infty$ and
$0<t_k\leq r_*^{-2}$.  This proves positivity of the two Hankel families displayed above.  In particular,
$\kappa_2=2S_1$, $\kappa_4=-12S_2$, and $\kappa_6=240S_3$, so
\begin{equation}\label{eq:G}
 G_a:=\kappa_2(a)\kappa_6(a)-\frac{10}{3}\kappa_4(a)^2
 =480(S_1S_3-S_2^2)>0.
\end{equation}
To verify the last sign, first expand finite partial sums and pair
the terms with indices $(i,j)$ and $(j,i)$; then pass to the limit:
\[
 S_1S_3-S_2^2=\sum_{i<j}t_it_j(t_i-t_j)^2.
\]
Absolute convergence justifies the expansion, and infinitely many
distinct zero heights ensure that at least one summand is strictly
positive.  Repeated zeros cause no difficulty, and simplicity is
not required.

The same cumulant combination, up to an overall factor of three, is displayed in Min and Kang \cite[Sec.~5.1]{MinKang}; we do not claim priority for it.

\subsection{A two-field matrix observation and its conditioning}
Fix $a=0$, $0<h_1<h_2$, and set $y_i=h_i^2$. The two-by-two matrix in
\eqref{eq:Loewner} can be measured without numerically differentiating a response curve:
\begin{equation}\label{eq:matrix-observed}
 L_{ii}=\frac{m_i-h_i\chi_i}{4h_i^3},\qquad
 L_{12}=\frac{m_1/(2h_1)-m_2/(2h_2)}{h_2^2-h_1^2},
 \quad m_i=\E_{q_{0,h_i}}X,\quad\chi_i=\Var_{q_{0,h_i}}X.
\end{equation}
Thus the first non-scalar matrix test requires two stationary ensembles and their means and variances. All parameters in this measurement prescription are well-defined without RH. Positivity is a necessary consequence of RH, but a successful two-field test is not a test of the universal quantifiers or proof of RH; indeed this particular matrix can also be checked directly from the known kernel, without assuming RH.

\paragraph{Fixed-design simultaneous inference.}
Choose fixed centers $c_i$ and radius $R>|c_i|$, and measure
$W_{ij}=\clip(X_{ij}-c_i,-R,R)$ from $N_i\geq2$ independent equilibrium draws per field.
Let $\widehat\mu_i=N_i^{-1}\sum_jW_{ij}$ and
$\widehat\nu_i=N_i^{-1}\sum_jW_{ij}^2$.
The physical targets are $m_i=c_i+\mu_i$ and
$\chi_i=\nu_i-\mu_i^2$, where $\mu_i=\E(X-c_i)$,
$\nu_i=\E(X-c_i)^2$ are \emph{unclipped}.
Write $s_{i,k}^2$ for the unbiased sample variance of $W^k$.
With $t=\log(16/\alpha)$, set
\begin{align}
 d_i&=\sqrt{2s_{i,1}^2t/N_i}+\frac{14Rt}{3(N_i-1)}+\beta_{i,1},\\
 e_i&=\sqrt{2s_{i,2}^2t/N_i}+\frac{7R^2t}{3(N_i-1)}+\beta_{i,2}.
\end{align}
Here $\beta_{i,1},\beta_{i,2}$ are deterministic clipping-bias bounds given below.
The empirical-Bernstein inequality and a union bound give, simultaneously with probability at least $1-\alpha$,
$|\mu_i-\widehat\mu_i|\leq d_i$ and
$|\nu_i-\widehat\nu_i|\leq e_i$ for both fields.
No independence between the mean and second-moment statistics is assumed.
Use $\widehat m_i=c_i+\widehat\mu_i$ and
$\widehat\chi_i=\widehat\nu_i-\widehat\mu_i^2$ in
\eqref{eq:matrix-observed} to define $\widehat L$. The same event implies
\begin{align}\label{eq:matrix-error}
 |L_{ii}-\widehat L_{ii}|&\leq E_{ii}:=
 \frac{d_i+h_i e_i+h_i(2|\widehat\mu_i|d_i+d_i^2)}{4h_i^3},\\
 |L_{12}-\widehat L_{12}|&\leq E_{12}:=
 \frac{d_1/(2h_1)+d_2/(2h_2)}{h_2^2-h_1^2}.
\end{align}
These are polynomial interval bounds, not an assumption that the plug-in variance or matrix estimator is unbiased.
A general sufficient sign rule is
$\lambda_{\min}(\widehat L)>\|E\|_2$ for positivity, or
$\lambda_{\min}(\widehat L)<-\|E\|_2$ for a negative eigenvalue.
For two-by-two matrices the entry rectangle is sharper. Define
\[
 d_i^-=\widehat L_{ii}-E_{ii},\quad d_i^+=\widehat L_{ii}+E_{ii},\quad
 b^+=|\widehat L_{12}|+E_{12},\quad b^-=\max(0,|\widehat L_{12}|-E_{12}),
\]
and $\Lambda(a,d,b)=(a+d-\sqrt{(a-d)^2+4b^2})/2$.
Monotonicity in the diagonal entries and in $|b|$ yields the simultaneous enclosure
\begin{equation}\label{eq:matrix-sign}
 \Lambda(d_1^-,d_2^-,b^+)\leq\lambda_{\min}(L)
 \leq\Lambda(d_1^+,d_2^+,b^-).
\end{equation}
A positive lower endpoint certifies this finite inequality; a negative upper endpoint would contradict RH only if the measured law were independently established to be the exact mathematical model, including all systematic error bounds.

\paragraph{Analytic clipping tails and numerical-law transfer.}
Use the appendix envelope $q_0(A+u)\leq D_Ae^{-\lambda_Au}$ for $u\geq0$.
Symmetry gives $M_0(h)=\E\cosh(hX)\geq1$.
For $A_\pm=R\pm c_i>0$, put
$\rho_\pm=\lambda_{A_\pm}\mp h_i$ and
$B_\pm=D_{A_\pm}e^{\pm h_iA_\pm}$.
If both rates are positive, integration of the two tilted tails gives
\begin{equation}
 \beta_{i,1}=\sum_{\pm}\frac{B_\pm}{\rho_\pm^2},\qquad
 \beta_{i,2}=\sum_{\pm}B_\pm
 \left(\frac{2R}{\rho_\pm^2}+\frac2{\rho_\pm^3}\right).
\end{equation}
The first bounds the absolute first-moment clipping bias and the second bounds the nonnegative second-moment clipping bias. Thus unbounded target moments are covered, not silently replaced by clipped ones.
If an independently justified numerical-law CDF discrepancy is at most $\varepsilon_i$, increase $d_i$ by $2R\varepsilon_i$ and $e_i$ by $2R^2\varepsilon_i$, the respective total-variation bounds for these clipped observables. Detector and model errors require additional independently established bounds; they are not covered by the sampling radius.

\paragraph{Conditioning and an actual numerical-law pilot.}
Deterministic kernel quadrature gives the following design comparison (rounding and quadrature convergence are numerical evidence, not interval arithmetic):
\begin{center}
\begin{tabular}{rrrr}
\toprule
$(h_1,h_2)$&$\lambda_{\min}(L)$&$\lambda_{\max}(L)$&condition number\\\midrule
$(5,8)$&$3.96298\,10^{-8}$&$5.54510\,10^{-5}$&$1399.23$\\
$(10,50)$&$4.00950\,10^{-7}$&$2.11380\,10^{-5}$&$52.72$\\\bottomrule
\end{tabular}
\end{center}
For the wider pair the predicted matrix is
\[
 L=\begin{pmatrix}2.03247255\,10^{-5}&4.02537557\,10^{-6}\\
 4.02537557\,10^{-6}&1.21423193\,10^{-6}\end{pmatrix}.
\]
A fixed-design pilot used $10^6$ inverse-CDF draws per field, seed $2026110307$, centers $(0.41,1.07)$, $R=1.5$, and $\alpha=0.05$.
Analytic rounded clipping bounds for the two fields are, respectively,
\[
 (\beta_1,\beta_2)\leq(6.596\,10^{-17},2.001\,10^{-16}),\quad
 (6.109\,10^{-14},1.853\,10^{-13}).
\]
The observed matrix was
\[
 \widehat L=\begin{pmatrix}2.04417887\,10^{-5}&4.02748860\,10^{-6}\\
 4.02748860\,10^{-6}&1.21563454\,10^{-6}\end{pmatrix}.
\]
Under the \emph{ideal exact-law sampling assumptions},
\eqref{eq:matrix-sign} gives a positive $95\%$ enclosure for the smallest eigenvalue.
The actual draw law was a floating-point inverse-CDF table, not exact $q_{0,h}$ sampling. A separate interval/Riemann audit of its CDF gives a conservative discrepancy envelope $2\cdot10^{-4}$ for each field, conditional on the stated interval elementary-function, floating-point interpolation, and ideal independent random-input assumptions. Transferring that envelope as described above gives the simultaneous 95\% enclosure
$[1.59748,6.17327]\,10^{-7}$.
This is a certificate for the finite mathematical matrix under those numerical and sampling assumptions, not a bit-level verification of the pseudorandom generator or a claim about a physical trap. The accompanying code, random seeds, checksums, and recorded output support reproduction; the full data archive additionally retains the raw draws. A crude operator-norm enclosure is inconclusive for the same data; the two-by-two rectangle calculation uses the markedly unequal entry-error scales.

\paragraph{Physical scale and limitations.}
At the physical units in \eqref{eq:units}, the wider design uses forces
$0.2057167$ and $1.0285835$ pN. Predicted mean positions are approximately
$81.5270$ and $214.4172$ nm, with standard deviations $36.13$ and $19.29$ nm.
The higher field probes a farther part of the target potential than the earlier weak-field examples: assuming that a calibrated device realizes $q_0$ there is an additional extrapolation, not demonstrated by the numerical pilot.
The stationary process has spectral gap at least $9$ in dimensionless time by the same curvature bound, uniformly in $h$. Hence for every square-integrable centered observable,
$0\leq\operatorname{Cov}(g(X_0),g(X_t))\leq e^{-9t}\operatorname{Var}(g)$,
and a continuous-time average over duration $T$ has variance at most
$2\operatorname{Var}(g)/(9T)$.
This aids planning but does not turn time-correlated observations into independent data or justify substituting an effective sample size into the displayed empirical-Bernstein theorem. A trajectory-based experiment needs its own dependence-aware concentration argument, burn-in analysis, detector calibration, and force-law error budget. The pilot is an equilibrium numerical-law simulation, not a Langevin trajectory or laboratory measurement.

\section{Unconditional zero-field signs and the low-order gap}
The theta-kernel moment inequalities of Csordas and Varga \cite{CsordasVarga} already provide unconditional Tur\'an inequalities for every real quadratic tilt; we do not identify those inequalities with the particular gap below without a separate derivation. The previous section must not be read as making every displayed low-order consequence RH-dependent. We record two stronger calibrations at $a=0$.

\begin{proposition}\label{prop:signs}
At $a=0$, all zero-field cumulants satisfy
\begin{equation}\label{eq:allsigns}
 (-1)^{j-1}\kappa_{2j}(0)>0\qquad(j\geq1)
\end{equation}
unconditionally. Moreover $G_0>0$ can be established directly by the interval integral certificate described below, without RH.
\end{proposition}
\begin{proof}
For each zeta zero $\rho=\beta+i\gamma$ with $\gamma>0$, put
$z_\rho=\gamma-i(\beta-1/2)$, counting multiplicities.
These are exactly the right-half-plane zeros of $\phi_0$.
The locally absolutely convergent product expansion from
Lemma~\ref{lem:analyticproduct} gives
\[
 \kappa_{2j}(0)=2(-1)^{j-1}(2j-1)!S_j,
 \qquad S_j=\sum_{\gamma>0}\Rea z_\rho^{-2j}.
\]
This formula does not require RH: the nonreal terms occur in
conjugate pairs, and the inverse-square summability also implies
summability of every higher even inverse power.

We first bound the total inverse-square weight without a Li
coefficient.  By Lemma~\ref{lem:zeros}, $\gamma>2$ and
$|\beta-1/2|\leq1/2$.  Set
$t=(\beta-1/2)^2/\gamma^2<1/16$.  Then
\[
 \Rea z_\rho^{-2}
 =\gamma^{-2}\frac{1-t}{(1+t)^2}
 \geq\frac{240}{289}\gamma^{-2}.
\]
Here $(1-t)/(1+t)^2$ is decreasing for $0\leq t\leq1/16$,
and its value at $1/16$ is $240/289$.
The variance identity gives $S_1=v_0/2$, while the curvature and
Brascamp--Lieb estimate used in Lemma~\ref{lem:zeros} gives
$v_0\leq1/18$.  Consequently
\[
 W:=\sum_{\gamma>0}\gamma^{-2}
 \leq\frac{289}{240}S_1\leq\frac{289}{8640}<0.05.
\]

The cosine integral in Table~\ref{tab:interval} is negative at
frequency $15$ and positive at zero.  Continuity therefore gives a
real Fourier zero $\gamma_*\in(0,15)$.
The finite-height theorem of Platt and Trudgian establishes that
all zeta zeros with $0<\gamma\leq T=3\cdot10^{12}$ lie on the
critical line \cite{PlattTrudgian}.  Their contribution to $S_j$
is positive and includes more than $15^{-2j}$ from $\gamma_*$.
For $j\geq2$ the remaining, potentially nonreal contribution obeys
\[
 \left|\sum_{\gamma>T}\Rea z_\rho^{-2j}\right|
 \leq\sum_{\gamma>T}|z_\rho|^{-2j}
 \leq\sum_{\gamma>T}\gamma^{-2j}
 \leq W T^{2-2j},
\]
because $|z_\rho|\geq\gamma$.
Its ratio to $15^{-2j}$ is strictly less than
\[
 \frac{0.05\cdot15^4}{T^2}
 \left(\frac{15}{T}\right)^{2j-4}
 \leq2.8125\cdot10^{-22}<1.
\]
Thus $S_j>0$ for every $j\geq2$.  The remaining case $j=1$
is $\kappa_2(0)=v_0>0$.
This argument uses the proved finite-height theorem and the
specified integral enclosure, rather than full RH.

For the gap, symmetry gives
$\kappa_2=\mu_2$,
$\kappa_4=\mu_4-3\mu_2^2$, and
$\kappa_6=\mu_6-15\mu_2\mu_4+30\mu_2^3$.
Substitution and cancellation yield exactly
\begin{equation}\label{eq:Gmoments}
 G_0=\mu_2\mu_6+5\mu_2^2\mu_4-\frac{10}{3}\mu_4^2.
\end{equation}
The interval integration described below and in
Appendix~\ref{app:integrals} encloses the positive mass and the
moments through degree six.  Normalization by the positive mass
interval and interval evaluation of this polynomial give the
strictly positive lower endpoint in Table~\ref{tab:interval}.
Accordingly $G_0>0$ follows under the explicitly stated
interval-arithmetic assumptions.  Ordinary floating-point values
alone are not used to establish that sign.
\end{proof}

\begin{table}[ht]
\centering\small
\begin{tabular}{p{0.46\linewidth}p{0.44\linewidth}}
\toprule
Quantity & Outward-enclosing interval\\\midrule
$\int_0^\infty K(x)\cos(15x)\dd x$ & $[-2.398514\cdot10^{-4},-1.129977\cdot10^{-4}]$\\
$G_0$ & $[2.937152\cdot10^{-7},1.578010\cdot10^{-6}]$\\
\bottomrule
\end{tabular}
\caption{Direct interval-box certificates. Rounded endpoints are widened from the machine-readable results. The computation assumes correct enclosure by \texttt{mpmath.iv} elementary functions; it is not a proof-assistant-checked calculation.}\label{tab:interval}
\end{table}

For completeness, the certificate integrates on $[0,2]$ using $40000$ closed boxes of width $1/20000$, evaluating the first three positive kernel terms by natural interval arithmetic at 40 decimal digits. Each box integral is enclosed by width times the interval range, not by an unproved quadrature remainder. The omitted series and spatial tails for moments through degree six have explicit positive majorants, detailed in Appendix~\ref{app:integrals}. Normalize moment intervals by the positive mass interval and evaluate \eqref{eq:Gmoments} by interval arithmetic. The cosine certificate uses the same boxes and the absolute unweighted tail majorant. The accompanying scripts retain all constants and output enclosures.

High-precision ordinary evaluation gives
\begin{align*}
 \kappa_2&=0.04620998623083794,\\
 \kappa_4&=-0.0004460711914232362,\\
 \kappa_6&=0.00003460174353623359,\\
 G_0/\kappa_2^4&=0.2052033682440877.
\end{align*}
These central values are numerical evaluations, distinct from the interval proof of positivity.

\paragraph{A simple logical counterexample.}
Let $P(X=1)=P(X=-1)=p/2$ and $P(X=0)=1-p$. Then
$M(h)=1-p+p\cosh h$ and
$G=p^2(15p-7)/3$. At $p=12/25$, $\kappa_4<0$, $\kappa_6>0$, and $G=48/3125>0$, but the zeros have nonzero real parts $\pm\operatorname{arcosh}(13/12)$. Thus initial alternating signs and a positive gap do not imply imaginary-axis zeros. At $p=9/20$, the same two signs hold while $G=-27/1600<0$. Even the gap is not a consequence of those signs alone.

\section{Reversible dynamics and physical normalization}
Consider the dimensionless diffusion
\begin{equation}\label{eq:SDE}
 \dd X_t=\{b_0(X_t)+aX_t+h/2\}\dd t+\dd W_t
 =\tfrac12(\log q_{a,h})'(X_t)\dd t+\dd W_t.
\end{equation}
It is nonexplosive, has invariant law $q_{a,h}$, and is reversible with generator
$\mathcal L=\tfrac12\partial_x^2+\tfrac12(\log q_{a,h})'\partial_x$.
Smooth, strongly inward double-exponential tails ensure the required existence and integrability even where the global curvature lower bound is not positive.

For a physical coordinate $r=\ell X$ in a thermal environment, let mobility be $\mu$, diffusivity $D=\mu k_BT$, and physical time be $t_{\rm p}=\tau t$ with
\begin{equation}\label{eq:units}
 \tau=\frac{\ell^2}{2D},\qquad
 U_{a,h}(r)=-k_BT\log q_0(r/\ell)
            -k_BT a(r/\ell)^2-k_BT h(r/\ell).
\end{equation}
Then \eqref{eq:SDE} is the overdamped equation
$\dd r=-\mu U_{a,h}'(r)\dd t_{\rm p}+\sqrt{2D}\dd W_{t_{\rm p}}$.
The applied constant force is $f=k_BT h/\ell$. A nonzero $a$ is a change in a quadratic potential, not a change in temperature. This exact potential is an idealized target for a calibrated trap or feedback device; no realization is asserted.

At $T=298\,$K and $\ell=0.2\,\mu$m, $h=5$ corresponds to $f=0.10285835\,$pN. Numerical equilibrium evaluation gives
\begin{center}
\begin{tabular}{rrrr}
\toprule
$h$ at $a=0$ & $m_0(h)$ & $\chi_0(h)$ & Physical mean and standard deviation\\\midrule
0 & 0 & 0.0462099862 & $0\,$nm, $42.9930\,$nm\\
5 & 0.2225654441 & 0.0414097686 & $44.5131\,$nm, $40.6988\,$nm\\
8 & 0.3389385310 & 0.0361062149 & $67.7877\,$nm, $38.0033\,$nm\\\bottomrule
\end{tabular}
\end{center}
Thus the mean shifts while equilibrium fluctuations narrow. These are predictions for \eqref{eq:units}, not measured bead data. The conversion uses the exact SI value of $k_B$; selecting a physical diffusivity sets the time unit independently.

\paragraph{Well-posedness and the stationary semigroup.}
We take each initial state to be independent of future Brownian increments.
Write $b_{a,h}=b_0+ax+h/2$. This drift is smooth and locally Lipschitz.
Its inward tail implies that, for suitable constants $C,c>0$ depending on
$a,h$, $2xb_{a,h}(x)\leq C-cx^2$. Applying the stopped It\^o formula
to $1+x^2$, then Gronwall and the usual exit-time estimate, proves
nonexplosion of the unique local strong solution. Applying the same
argument to $1+x^{2p}$ proves finite moments on bounded time intervals
for every integer $p\geq1$ when the initial law has these moments.
Moreover $q_{a,h}'=2b_{a,h}q_{a,h}$, and the tails make all boundary terms
below vanish. For smooth compactly supported $u,v$,
\[
 \int u\mathcal Lv\,q_{a,h}\,\dd x
 =-\frac12\int u'v' q_{a,h}\,\dd x
 =\int v\mathcal Lu\,q_{a,h}\,\dd x.
\]
Thus the diffusion semigroup is the conservative symmetric semigroup
associated with this Dirichlet form, and $q_{a,h}$ is invariant.
The rapidly decaying density also puts the centered coordinate in the
form domain and gives its Dirichlet energy $1/2$.

\begin{theorem}\label{thm:dynamics}
Initialize \eqref{eq:SDE} in $q_{a,h}$, with $a\geq0$, and let
$C_h(t)=\Cov(X_0,X_t)$ and $S_h(t)=\chi_a(h)-C_h(t)$. Then
\begin{equation}\label{eq:dynamicbound}
 0\leq S_h(t)\leq\chi_a(h)
 \left(1-e^{-t/(2\chi_a(h))}\right)
 \leq v_a\left(1-e^{-t/(2v_a)}\right).
\end{equation}
Here $S_h(t)$ is the derivative of the mean at time $t$ under an infinitesimal field step from $h$ to $h+\varepsilon$, with the initial law kept at $q_{a,h}$.
If a deterministic $u\in L^2[0,T]$ adds drift $u(t)/2$ while keeping the same initial law, then
\begin{equation}\label{eq:energy}
 \left|\E^uX_T-m_a(h)\right|
 \leq\frac{\sqrt{v_a}}2\|u\|_{L^2[0,T]}.
\end{equation}
Both conclusions are unconditional.
\end{theorem}
\begin{proof}
Put $q=q_{a,h}$ and $f(x)=x-m_a(h)$. Reversibility and stationarity give
$C_h(t)=\langle f,P_tf\rangle_{L^2(q)}$. The spectral theorem for the
nonnegative self-adjoint operator $-\mathcal L$ therefore supplies a
positive finite measure $\nu_f$ on $[0,\infty)$ such that
\[
 C_h(t)=\int e^{-\lambda t}\,\dd\nu_f(\lambda),\qquad
 \nu_f([0,\infty))=\|f\|_2^2=\chi_a(h),\qquad
 \int\lambda\,\dd\nu_f(\lambda)=\frac12\int(f')^2q\,\dd x=\frac12.
\]
The last identity is a form-domain identity and does not require assuming
that the coordinate is bounded. Since $\chi_a(h)>0$, Jensen applied to
the probability measure $\nu_f/\chi_a(h)$ yields
\[
 C_h(t)\geq\chi_a(h)\exp\{-t/(2\chi_a(h))\}.
\]
The same positive-measure representation gives $0\leq C_h(t)\leq\chi_a(h)$.
Consequently the first two inequalities in \eqref{eq:dynamicbound} hold.
For $v>0$, differentiation gives
\[
 \frac{\dd}{\dd v}\{v(1-e^{-t/(2v)})\}
 =1-\left(1+\frac{t}{2v}\right)e^{-t/(2v)}\geq0,
\]
because $e^x\geq1+x$. The last inequality follows from
$\chi_a(h)\leq v_a$.

For clarity, the field step changes the generator to
$\mathcal L_\varepsilon=\mathcal L+(\varepsilon/2)\partial_x$ but leaves
its initial distribution equal to $q$. Duhamel differentiation,
stationarity of the unperturbed semigroup, and the Dirichlet identity
with $f'=1$ give
\begin{align*}
 \left.\partial_\varepsilon\E^{h+\varepsilon}X_t\right|_{\varepsilon=0}
 &=\frac12\int_0^t\int(P_sf)'q\,\dd x\,\dd s\\
 &=-\int_0^t\langle f,\mathcal LP_sf\rangle_{L^2(q)}\,\dd s\\
 &=-\int_0^t C_h'(s)\,\dd s=C_h(0)-C_h(t).
\end{align*}
One first applies these identities to smooth truncations and then removes
the truncation using the moment bounds and the Dirichlet-form identity.
Equivalently, differentiability of the field-step law follows from
Girsanov with the deterministic drift increment $\varepsilon/2$.

For the control statement, let $P^0,P^u$ be the path laws on $[0,T]$
with the same initial law $q$. The Girsanov shift is $u/2$, so its
Novikov integral is deterministic and finite. Under $P^u$ the log
likelihood ratio has expectation
\[
 D(P^u\Vert P^0)=\frac12\int_0^T\left(\frac{u(s)}2\right)^2\dd s
 =\frac18\|u\|_2^2.
\]
For every real $r$, Taylor's integral formula applied to
$\log M_a(h+r)-\log M_a(h)-rm_a(h)$ gives
\[
 \log\E_q e^{r(X-m_a(h))}
 =r^2\int_0^1(1-s)\chi_a(h+sr)\,\dd s
 \leq\frac{v_a r^2}{2}.
\]
If $Q$ is the controlled terminal law, the entropy variational inequality
and data processing imply, for $r>0$,
\[
 \pm\{\E_QX-m_a(h)\}
 \leq\frac{D(Q\Vert q)}r+\frac{v_ar}{2}
 \leq\frac{D(P^u\Vert P^0)}r+\frac{v_ar}{2}.
\]
Optimizing in $r$ proves \eqref{eq:energy}; the zero-entropy case follows
by taking $r\downarrow0$. The quantity used here is path relative
entropy, not heat or thermodynamic work.
\end{proof}

\section{Route C quadratic deformation and a variance quench}
The parameter $a$ has an exact transform interpretation:
\begin{equation}\label{eq:heat}
 \partial_a\phi_a(z)=-\partial_z^2\phi_a(z)-v_a\phi_a(z),
 \qquad \frac{\dd v_a}{\dd a}=\Var_{q_a}(X^2)>0.
\end{equation}
The first equation follows by differentiation and normalization; the second is the exponential-family covariance identity. Neither uses RH. It is essential to distinguish this deformation parameter from elapsed diffusion time. In the conventional heat notation its value is $4a$, as in \eqref{eq:scaling}.

\begin{theorem}\label{thm:quench}
Fix $0<A<9$, put $c=9-A$, and start
\begin{equation}\label{eq:quench}
 \dd X_t=\{b_0(X_t)+AX_t\}\dd t+\dd W_t,
 \qquad X_0\sim q_0.
\end{equation}
Let $V(t)=\E X_t^2$ and $v_A=\E_{q_A}X^2$. Then $V(t)$ is nondecreasing,
\begin{equation}\label{eq:quenchbound}
 0\leq v_A-V(t)\leq(v_A-v_0)e^{-2ct},
 \qquad V'(0)=2Av_0.
\end{equation}
\end{theorem}
\begin{proof}
The curvature bound implies $b_A'(x)\leq-c$ and oddness gives $b_A(0)=0$.
Thus $b_A(r)\leq-cr$ for $r\geq0$. Tanaka's formula shows that the
absolute-value process has, in law, the reflected-diffusion equation
\[
 \dd R_t=b_A(R_t)\dd t+\dd B_t+\dd L_t,\qquad R_t\geq0,
 \qquad \int_0^t\mathbf1_{\{R_s>0\}}\dd L_s=0.
\]
Here $B$ is a Brownian motion; it need not be the Brownian motion driving
the original signed coordinate. Locally Lipschitz reflected equations
have pathwise comparison under a common reflected driving noise:
ordered starting points remain ordered when the drifts are ordered.
To see the comparison, apply the positive-part formula to the difference
of the two solutions stopped in a compact interval. Their Brownian parts
cancel. At a first crossing the ordered-drift condition prevents a
positive crossing, and reflection of the larger process cannot occur
while its position is strictly larger. Local Lipschitz continuity and
Gronwall finish the stopped comparison, after which nonexplosion removes
the stops.

Let $\nu_a$ be the law of $|X|$ under $q_a$. Its likelihood ratio satisfies
$\dd\nu_A/\dd\nu_0=Ce^{Ar^2}$, an increasing function of $r$. Hence
$\nu_0\leq_{\rm st}\nu_A$. For example, for bounded increasing $f$, the
inequality follows directly from the nonnegative covariance of the two
increasing functions $f(r)$ and $e^{Ar^2}$ under $\nu_0$.
As $b_A(r)\geq b_0(r)$ for $r\geq0$, comparison with a stationary
reflected $b_0$ process yields $\nu_0P_A^s\geq_{\rm st}\nu_0$.
The same-drift comparison makes $P_A^t$ order preserving; applying it
to this inequality gives
\[
 \nu_0P_A^{s+t}=(\nu_0P_A^s)P_A^t\geq_{\rm st}\nu_0P_A^t.
\]
Likewise $\nu_0\leq_{\rm st}\nu_A$ and invariance of $\nu_A$ imply
$\nu_0P_A^t\leq_{\rm st}\nu_A$. Approximating the increasing function
$r^2$ by bounded increasing functions proves that $V(t)$ is
nondecreasing and $V(t)\leq v_A$.

For the rate, use the quantile coupling of initial radii
$R_0\sim\nu_A$, $S_0\sim\nu_0$ with $R_0\geq S_0$, and thereafter
use the same reflected noise and drift $b_A$. Then $R_t\geq S_t$ and
$R_t\sim\nu_A$. In It\^o's formula for $r^2$ the reflection term
$2R_t\dd L_t$ vanishes, since $L$ increases only at zero. Its generator
value is $J(r)=1+2rb_A(r)$ and
\[
 J'(r)=2b_A(r)+2rb_A'(r)\leq-4cr.
\]
Integration from $S$ to $R$ gives
$J(R)-J(S)\leq-2c(R^2-S^2)$. Thus, writing
$d(t)=\E(R_t^2-S_t^2)=v_A-V(t)$, the stopped It\^o formula and then
removal of stops give $d'(t)\leq-2cd(t)$ almost everywhere. Gronwall
proves $d(t)\leq d(0)e^{-2ct}$. Polynomial moments are controlled by
\[
 \mathcal L_A x^{2p}
 \leq-2pcx^{2p}+p(2p-1)x^{2p-2},
\]
which justifies the moment limits. Integrability of the drift terms
also follows from the stationary upper process and the monotonicity
of $J$ on $[0,\infty)$.
Finally, stationarity at $a=0$ gives
$\E_{q_0}[1+2Xb_0(X)]=0$. Therefore
\[
 V'(0)=\E_{q_0}[1+2Xb_A(X)]
 =\E_{q_0}[1+2Xb_0(X)]+2A\E_{q_0}X^2=2Av_0.
\]
\end{proof}

At $A=2$, ordinary quadrature and two converging PDE discretizations give
\[
 v_0=0.04620998623,\quad
 V(0.1)\approx0.05405111824,\quad
 v_2=0.05519528393.
\]
The predicted variance increase by $t=0.1$ is about $16.97\%$. With $\ell=0.2\,\mu$m, the root-mean-square coordinate increases from $42.9930\,$nm to approximately $46.4978\,$nm, approaching $46.9874\,$nm. For illustration, choosing $D=0.4905\,\mu\mathrm m^2/\mathrm s$ gives $\tau\approx40.77\,$ms and $t=0.1$ corresponds to about $4.08\,$ms. This chosen diffusivity is a model parameter, not a claimed measurement. Agreement between finite-volume and transformed-operator discretizations is convergence evidence, not an interval-certified PDE solution.

\paragraph{A matched harmonic comparator.}
As a specified reference, start from the same $q_0$ law but replace $b_0$ by the harmonic drift $-\gamma_0 x$, where $\gamma_0=(2v_0)^{-1}$, and apply the identical additional drift $Ax$. Its variance is explicit:
\begin{equation}\label{eq:oucomp}
 V_{\rm OU}(t)=\frac1{2(\gamma_0-A)}+\left(v_0-\frac1{2(\gamma_0-A)}\right)e^{-2(\gamma_0-A)t},\qquad A<\gamma_0.
\end{equation}
For $A=2,t=0.1$, this yields $V_{\rm OU}=0.05489276589$, versus $V=0.05405111824$ for the theta model, a difference of $0.00084165$. Matching the initial variance and quench coefficient does not match the full potential or its post-quench equilibrium variance. This is a controlled model comparator, not a literature benchmark or an experimental observation. Two independent discretization formulations agree after extrapolation; this is numerical convergence evidence, not an interval error certificate.
\begin{figure}[htbp]
\centering
\includegraphics[width=\linewidth]{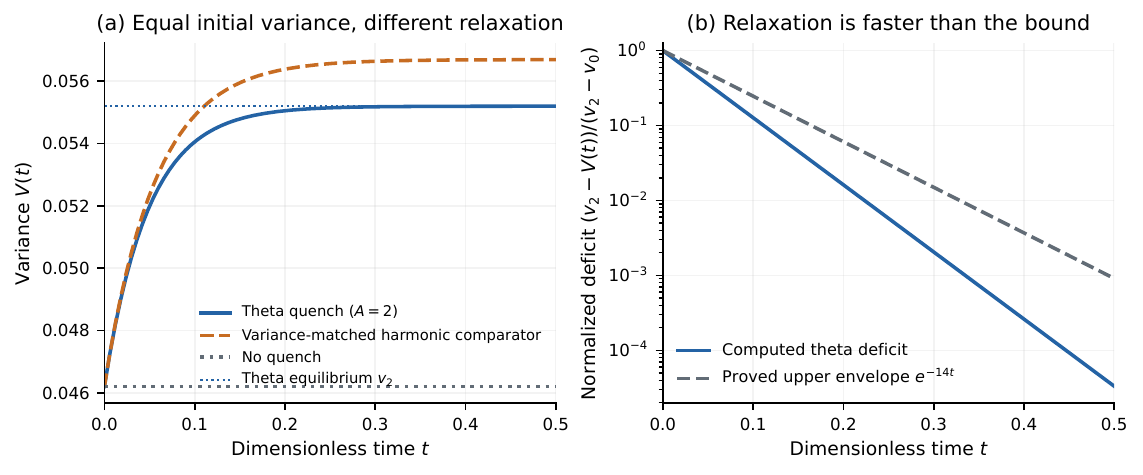}
\caption{Numerical predictions for the ideal quench and variance-matched harmonic comparator, starting from the same initial law. Left: variance response. Right: normalized theta-model variance deficit and the proved bound $e^{-14t}$. Curves use deterministic numerical discretization, not measured data.}\label{fig:quench}
\end{figure}

\subsection{A finite implementation error bound}
Finite feedback changes the model. Timing and observation corrections for virtual feedback potentials were analyzed by Jun and Bechhoefer \cite{JunBechhoefer}; the following bound addresses a specified clipped nonlinear implementation rather than replacing that literature. Fix $0<A<9$, $\Delta>0$, $B>0$, $\sigma\geq0$, and an integer delay index $m\geq0$ (distinct from the equilibrium mean $m_a(h)$). For example clip the \emph{total} drift to $g=\clip(b_A,-B,B)$ and apply delayed noisy sample-and-hold actuation,
\begin{equation}\label{eq:implemented}
 \dd Z_t=g(Z_{(k-m)\Delta}+\sigma\xi_{k-m})\dd t+\dd W_t,
 \quad k\Delta\leq t<(k+1)\Delta.
\end{equation}
For negative indices set the state to $Z_0$ and initialize independent measurement noise. Each measurement is retained until actuation, not freshly remeasured then. Let $\xi_j$ be independent standard normals, independent of the initial state and Brownian motion. Couple \eqref{eq:implemented} with the ideal process \eqref{eq:quench} using the same Brownian noise. Allow imperfect preparation, with $e_0=\|Z_0-X_0\|_2$, where $X_0\sim q_0$. Set
\[
 K_B=\operatorname{Lip}(g),\quad
 \epsilon_B(t)=\|(|b_A(X_t)|-B)_+\|_2,\quad
 h_*=(m+1)\Delta.
\]
Since $g$ is nonincreasing, the common-noise difference and Cauchy--Schwarz imply
\begin{equation}\label{eq:implementationbound}
 e(t):=\|Z_t-X_t\|_2
 \leq e_0+tK_B(Bh_*+\sqrt{h_*}+\sigma)
       +\int_0^t\epsilon_B(s)\dd s.
\end{equation}
Indeed decompose the drift error into $g(Z_t)-g(X_t)$, the held-input error, and $g(X_t)-b_A(X_t)$. The first has nonpositive inner product with $Z_t-X_t$; the delayed state increment has $L^2$ norm at most $Bh_*+\sqrt{h_*}$, and the noise contribution is at most $\sigma$. Integrating the upper derivative of the $L^2$ norm proves the estimate. No independence between already actuated measurement noise and the later state is assumed.

The ideal mean is zero by symmetry. For second moments,
\begin{equation}\label{eq:varianceerror}
 |\E Z_t^2-\E X_t^2|
 \leq e(t)\{2\sqrt{\E X_t^2}+e(t)\}.
\end{equation}
The same right-hand side bounds $|\Var Z_t-\Var X_t|$: centering is an $L^2$ contraction and the standard deviation is the $L^2$ distance to constants. Thus no symmetry of the implemented initial law is needed. This controls the unbounded observable directly, rather than incorrectly inferring moment accuracy from total variation alone. The ideal second moment obeys
$\E X_t^2\leq e^{-2ct}v_0+(1-e^{-2ct})/(2c)$.
Finally $q_0=Z(A,0)e^{-Ax^2}q_A\leq Z(A,0)q_A$ and positivity of the Markov semigroup give
\begin{equation}\label{eq:clippingtail}
 \sup_{t\geq0}\epsilon_B(t)^2
 \leq\int_\R q_0(x)e^{Ax^2}(|b_A(x)|-B)_+^2\dd x.
\end{equation}
This is an exact one-dimensional tail bound. The observation equation above is a point-sampling model, not a camera-exposure model. If instead the noiseless observation at $j\Delta$ is the backward exposure average $\bar Z_j=\epsilon^{-1}\int_{j\Delta-\epsilon}^{j\Delta}Z_s\,\mathrm ds$, and the process has the same unit Brownian diffusion and drift bounded by $B$ over the required prehistory, the bounded drift gives
\begin{equation}\label{eq:exposure}
 \|\bar Z_j-Z_{j\Delta}\|_2\leq B\epsilon/2+\sqrt{\epsilon/3}.
\end{equation}
Indeed the drift part is at most $B\epsilon/2$ pathwise, while the average Brownian increment has variance $\epsilon/3$. Thus the right side of \eqref{eq:implementationbound} remains valid after adding $tK_B(B\epsilon/2+\sqrt{\epsilon/3})$; the measurement noise term remains separately included. This bound allows overlap of exposures and temporal correlation, but assumes the stated additive observation model and bounded actuation. It does not include hydrodynamic memory, calibration error, or exposure-dependent detector noise. Clipped drift is not globally $c$-strongly dissipative, so the ideal contraction constant must not be assigned to it. At $A=2,B=30,\Delta=10^{-3},m=1,\sigma=0.01$, the elementary uniform bound is too conservative to certify a small error relative to the quench signal. Numerical simulations of this implementation and a rigorous small-error certificate are different claims. No stationary distribution for the delayed noisy process is asserted.

\paragraph{Details of the elementary implementation bound.}
Let $D_t=Z_t-X_t$ and let $j(t)=\max\{0,(k-m)\Delta\}$ on the $k$th
hold interval, retaining the original negative-index noise label when
necessary. Common noise cancels in $\dd D_t$. Add and subtract $g(Z_t)$
and $g(X_t)$. The contribution
$D_t[g(Z_t)-g(X_t)]$ is nonpositive. With
$e(t)=\|D_t\|_2$, Cauchy--Schwarz therefore gives, in the upper derivative
sense (or after regularizing the square root),
\[
 D^+e(t)\leq K_B(\|Z_t-Z_{j(t)}\|_2+\sigma)+\epsilon_B(t).
\]
Bounded actuation and Minkowski give
\[
 \|Z_t-Z_{j(t)}\|_2
 \leq\int_{j(t)}^t\|g(\hbox{held input})\|_2\,\dd s
       +\|W_t-W_{j(t)}\|_2
 \leq B h_*+\sqrt{h_*}.
\]
Integrating proves \eqref{eq:implementationbound}. Only the marginal
identity $\|\xi_j\|_2=1$ is used here; already actuated noises need not
be independent of the later state.
For any square-integrable coupled $X,Z$,
\[
 |\E Z^2-\E X^2|
 \leq\|Z-X\|_2\|Z+X\|_2
 \leq e(2\|X\|_2+e).
\]
Also centering is an orthogonal projection in $L^2$, so
$|\sqrt{\Var Z}-\sqrt{\Var X}|\leq e$. Hence
$|\Var Z-\Var X|\leq e(2\sqrt{\Var X}+e)$, which proves the stated
variance transfer without assuming symmetry of the implemented initial
law.

The domination used for the source terms is a measure inequality.
Indeed $q_0\,\dd x\leq Z(A,0)q_A\,\dd x$ and, for every bounded
nonnegative $f$, positivity and invariance yield
\[
 \E f(X_t)=\int P_A^tf\,q_0\,\dd x
 \leq Z(A,0)\int P_A^tf\,q_A\,\dd x
 =\int f(x)q_0(x)e^{Ax^2}\,\dd x.
\]
Monotone convergence extends this to the nonnegative squared source
functions in the implementation estimates.

For the exposure statement, the additional drift contribution to the
backward average is at most $B\epsilon/2$. Stochastic Fubini writes its
Brownian contribution, up to sign, as
\[
 \frac1\epsilon\int_{t-\epsilon}^t
       (s-t+\epsilon)\,\dd W_s,
\]
whose variance is
$\epsilon^{-2}\int_0^\epsilon u^2\dd u=\epsilon/3$.
This proves \eqref{eq:exposure}. A Brownian, bounded-drift prehistory is
required whenever an exposure window reaches before initialization;
the constant negative-index convention by itself does not supply that
prehistory. This elementary exposure modification does not establish
a finite-exposure version of the sharper grid-time delay proposition.

\paragraph{An analytically bounded fast-feedback design.}
In contrast to the preceding vacuous example, the same clipped-total-drift model admits a deliberately extreme theoretical design. At $A=2$, $B=30$, $m=1$, $\Delta=10^{-8}$, $\sigma=10^{-5}$ and exact shared preparation $e_0=0$, the analytic tail and Lipschitz estimates in Appendix~\ref{app:fast} give $K_B<80$ and $\sup_t\epsilon_B(t)<3.4\cdot10^{-5}$. Equation~\eqref{eq:implementationbound} then yields, at $t=0.1$, $e(t)<0.001220$ and
\begin{equation}\label{eq:fastbudget}
 |\Var Z_t-V(t)|<0.000654.
\end{equation}
Here we used $V(t)\leq1/14$. The exact quench theorem gives $V(0.1)-v_0\geq(v_2-v_0)(1-e^{-1.4})$; ordinary quadrature gives this lower expression approximately $0.00677$, whereas $V(0.1)-v_0\simeq0.00784$. The displayed numeric budget is obtained from elementary analytic inequalities with conservatively rounded ordinary arithmetic, not a machine-directed interval certificate for every decimal. The update period corresponds to about $0.408$ ns for the illustrative physical parameters, so this example proves a mathematical small-step regime rather than practical device feasibility. Any preparation error $e_0$, detector calibration error, exposure contribution, or mismatch in the assumed force law must be added to the relevant budget. It does not resolve the smaller difference from the matched harmonic comparator.

\subsection{One-step-delay cancellation for bounded total drift}
The very small update period above is sufficient but not necessary for a non-vacuous bound. For one-step delay one can retain cancellation in the discrete dynamics rather than summing all local Brownian errors in absolute value. Let $A=2$, $b=b_2$, $g=\clip(b,-B,B)$ and $G(x)=\E_\xi g(x+\sigma\xi)$, with $\xi\sim N(0,1)$. Write $\rho(x)=q_0(x)e^{2x^2}$ and define
\[
 C_1^2=\int\rho (b')^2,\quad C_0^2=\int\rho (bb'+b''/2)^2,\quad
 C_s^2=\int\rho(G-b)^2,
\]
where all integrals are over $\mathbb R$. Set $K=\operatorname{Lip}(g)$.

\begin{proposition}\label{prop:delay}
For the one-step delayed implementation \eqref{eq:implemented} with $m=1$, take $Z_{-\Delta}=Z_0=X_0\sim q_0$, define the formal $X_{-\Delta}=X_0$, and use independent stored Gaussian measurements, and a common Brownian path. If $N\geq1$ is an integer, $\alpha>0$, and $\Delta K\leq2\alpha/(\alpha^2+\alpha+1)$, then at $t_N=N\Delta$,
\begin{align}\label{eq:delaybound}
 \|Z_{t_N}-X_{t_N}\|_2
 &\leq\sqrt{1+\alpha}\bigg[\left(\frac{7N-6}{3}C_1^2\Delta^3+N\Delta^2K^2\sigma^2\right)^{1/2}\\
 &\hspace{15mm}+\frac{3N-2}{2}\Delta^2 C_0+N\Delta C_s\bigg].\nonumber
\end{align}
This result is unconditional on RH. It assumes instantaneous point observations, not exposure averaging, and exactly prepared $q_0$.
\end{proposition}
\begin{proof} The full calculation, including the adapted filtrations and the surviving covariance term, is given in Appendix~\ref{app:delayproof}.\end{proof}

\begin{figure}[htbp]
\centering
\includegraphics[width=\linewidth]{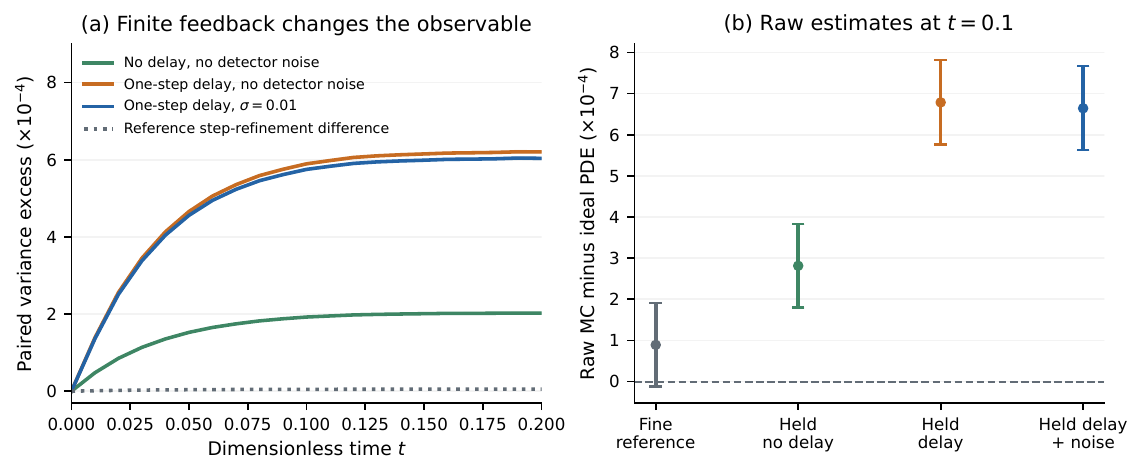}
\caption{Finite-feedback simulation. Left: paired excess in second moment relative to the fine Euler reference; the dotted curve diagnoses a refinement difference. Right: raw second-moment estimates relative to the deterministic ideal PDE, with 95\% Monte Carlo bars only. Curves assume the stated clipped feedback and do not represent measured device data.}\label{fig:feedback}
\end{figure}

For $B=30$, $\Delta=10^{-3}$, $\sigma=0.01$, $N=100$, and $\alpha=0.04$, the bounds in Appendix~\ref{app:delayconstants} give
\[
 K<70,\quad C_1<10.1,\quad C_0<26.7,\quad C_s<0.00135,
 \quad v_0<0.046330,\quad 0.055045<v_2<0.055346.
\]
The step-size criterion holds. Equation~\eqref{eq:delaybound} gives $\|Z_{0.1}-X_{0.1}\|_2<0.0128840$, hence $|\Var Z_{0.1}-V(0.1)|<0.006229$. Combining the quench theorem with the moment enclosures gives $V(0.1)-v_0>0.0065659$, so
\begin{equation}\label{eq:delaynumeric}
 \Var Z_{0.1}-v_0>0.000336.
\end{equation}
This is a non-vacuous, but loose, bound for the \emph{total-clipped} model at the same update, delay, and noise parameters as Fig.~\ref{fig:feedback}. It does not certify the much smaller numerical error or resolve the theta--harmonic difference. The enclosures depend on the supplied interval-function and floating aggregation assumptions; no physical calibration, finite exposure, or nonideal initialization is included.

\subsection{Passive confinement with bounded residual feedback}
A different architecture yields a more informative finite-resource certificate. The harmonic confinement is assumed to act continuously at the instantaneous position, while only a bounded residual is delayed and held. Let $0<A<9$, $c=9-A$, $r(x)=b_0(x)+9x$, and $r_R(x)=r(\clip(x,-R,R))$. The curvature estimate makes $r$ nonincreasing. With $B=\sup|r_R|$, $K=\operatorname{Lip}(r_R)$, define
\begin{equation}\label{eq:passive}
 \dd Z_t=[-cZ_t+r_R(Z_{(k-m)\Delta}+\sigma\xi_{k-m})]\dd t+\dd W_t,
 \qquad k\Delta\leq t<(k+1)\Delta.
\end{equation}
The negative-index state convention and stored independent Gaussian observation errors are as in \eqref{eq:implemented}. This is \emph{not} a total-force cap: the passive harmonic force is unbounded on the real line. Its accuracy, spatial range, detector exposure, and physical realization require separate justification.

\begin{proposition}\label{prop:passive}
Start the ideal quench and \eqref{eq:passive} with the same $X_0=Z_0\sim q_0$ and the same Brownian motion. Let $h_*=(m+1)\Delta$ and $Kh_*<1$. If
\[
 \epsilon^2\geq\int_{\mathbb R}q_0(x)e^{Ax^2}|r(x)-r_R(x)|^2\dd x,
\]
set
\[
 F=K[(c\sqrt{v_A}+B)h_*+\sqrt{h_*}+\sigma]+\epsilon,\quad
 \lambda=c(1-Kh_*),\quad E(t)=\frac{F(1-e^{-\lambda t})}{\lambda}.
\]
Then $\|Z_t-X_t\|_2\leq E(t)$ and, unconditionally on RH,
\begin{equation}\label{eq:passivebound}
 \E Z_t^2-v_0\geq(v_A-v_0)(1-e^{-2ct})-E(t)[2\sqrt{v_A}+E(t)].
\end{equation}
By symmetry $\E Z_t=0$, so the left side is a variance increase. The assertion does not imply monotonicity of the implemented variance.
\end{proposition}
\begin{proof}
Write $g(x)=-cx+r_R(x)$. Since $r_R$ is nonincreasing,
\[
 (z-x)[g(z)-g(x)]\leq-c(z-x)^2.
\]
On each hold interval the implemented equation is a linear equation
with bounded forcing, so it has a unique global strong solution with
finite moments. The ideal law has $\|X_t\|_2\leq\sqrt{v_A}$ by the
quench theorem, and measure domination gives
$\|r(X_t)-r_R(X_t)\|_2\leq\epsilon$.
For $j=\max\{0,(k-m)\Delta\}$, decompose the common-noise drift error as
\[
 g(Z_t)-g(X_t)
 +r_R(Z_j+\sigma\xi_{k-m})-r_R(Z_t)
 +r_R(X_t)-r(X_t).
\]
Dissipativity and Cauchy--Schwarz imply
\[
 D^+e(t)\leq-ce(t)+K\{\|Z_t-Z_j\|_2+\sigma\}+\epsilon.
\]
Since $t-j\leq h_*$, Minkowski and the bounded residual give
\begin{align*}
 \|Z_t-Z_j\|_2
 &\leq\int_j^t\{c\|Z_s\|_2+B\}\,\dd s+\sqrt{t-j}\\
 &\leq(c\sqrt{v_A}+B)h_*+ch_*\sup_{0\leq s\leq t}e(s)+\sqrt{h_*}.
\end{align*}
Thus $D^+e\leq-ce+cKh_*\sup_{s\leq t}e(s)+F$.
Here $e(0)=0$ and $\lambda=c(1-Kh_*)>0$.
The candidate $E(t)=F(1-e^{-\lambda t})/\lambda$ is increasing and
satisfies $E'=-\lambda E+F$. To make the comparison strict, fix
$\delta>0$ and use $E_\delta(t)=E(t)+\delta e^t$.
It starts strictly above $e$, is increasing, and satisfies
$E_\delta'>-\lambda E_\delta+F$.
At a first crossing, all earlier errors are at most
$E_\delta(t)$, so the error differential inequality prevents that
crossing. Equivalently the same argument can be applied after
regularizing the $L^2$ norm. Hence $e\leq E_\delta$; letting
$\delta\downarrow0$ gives $e\leq E$.
The second-moment transfer proved above now gives
\[
 \E Z_t^2\geq V(t)-E(t)(2\sqrt{v_A}+E(t)).
\]
Combining this with
$V(t)-v_0\geq(v_A-v_0)(1-e^{-2ct})$ proves
\eqref{eq:passivebound}.
Finally $q_0$ is even, $r_R$ is odd, and the joint driving law is invariant
under $(X_0,W,\xi)\mapsto(-X_0,-W,-\xi)$. Pathwise uniqueness implies
$Z_t\stackrel{\rm law}=-Z_t$, so its finite first moment is zero.
No monotonicity of the implemented variance, and no independence of an
already actuated noise from a later state, was used.
\end{proof}

For $A=2$, $R=1$, one may take $c=7$, $B=12$, $K=38$, and $\epsilon=0.002$; Appendix~\ref{app:passive} gives analytic bounds. At $\Delta=2\cdot10^{-6}$, $m=1$, $\sigma=0.001$, $t=0.1$, an interval-box calculation gives $v_0\in[0.04616354,0.04625647]$ and $v_2\in[0.05513614,0.05525449]$ with widened endpoints. It follows that $E(0.1)<0.008492$, the variance error is below $0.004065$, and
\begin{equation}\label{eq:passivenumber}
 \Var Z_{0.1}-v_0>0.00262.
\end{equation}
These inequalities assume correct enclosures by the specified interval library and the analytic tail bounds, not a verified hardware experiment. For illustration only, $\ell=10\,\mu$m and $D=0.5\,\mu$m$^2$/s give $\tau=100$ s, a $0.2$ ms update, $0.2$ ms one-step delay, $10$ nm input noise and $t=0.1$ equal to 10 s. These parameter choices do not establish a feasible or calibrated device. The result demonstrates a non-vacuous mathematical control budget, not an RH test or a guarantee at finite camera exposure.

\paragraph{Simulated feedback at a finite update rate.}
To gauge the size of systematic distortion, simulate the total-clipped implementation at $A=2$, $B=30$, update $\Delta=10^{-3}$, one update of delay, and noise standard deviation $\sigma=0.01$. A paired common-noise calculation with $2{,}097{,}152$ paths, 64 batches, and seed 2026110317 gives at $t=0.1$ a difference in second moments between the held implementation and a fine Euler reference of $0.000574983\pm0.000001864$ (nominal 95\% Monte Carlo error only). The fine reference uses step $2.5\cdot10^{-5}$. Comparing it with step $5\cdot10^{-5}$ changes the second moment by approximately $4.92\cdot10^{-6}$; neither the Monte Carlo interval nor the step comparison certifies residual discretization or lookup error. The implementation is a numerical model, not a physical realization. The numerical implementation shift is of the same order as the $0.00084165$ ideal theta-versus-harmonic difference, underscoring the need to calibrate feedback before attributing such a difference to the potential.

\section{Route A a tilted bridge response and its discretization}
At $a=0$ the mean has the direct representation
\begin{equation}\label{eq:A}
 m_0(h)=\frac{\E_P[Y^{h+1/2}\log Y]}{\E_PY^{h+1/2}},
 \qquad A(h):=\tfrac12m_0(h).
\end{equation}
For $h>0$, Theorem~\ref{thm:static} gives
$0<A(h)<hv_0/2$ and a strict susceptibility comparison. These statements are unconditional.

At $h=1$, ordinary high-precision evaluation gives
$A(1)=0.02306796403023129$.

A finite-grid bridge range is not the continuum range. In a computational audit, one million bridge paths were evaluated on nested meshes; the ordinary weighted-ratio estimates are shown in Table~\ref{tab:A}. The Monte Carlo standard error concerns the finite-grid estimand, not the missing extrema. At the finest listed mesh the discrepancy from the continuum value is about $39.5$ Monte Carlo standard errors. More independent paths at a fixed mesh do not remove this bias.
\begin{table}[ht]
\centering\small
\begin{tabular}{rrrr}
\toprule
Mesh intervals & Estimate of $A(1)$ & Monte Carlo SE & Estimate minus continuum\\\midrule
128 & $-0.015596623$ & $0.000127796$ & $-0.038664587$\\
512 & $0.004074491$ & $0.000122265$ & $-0.018993473$\\
2048 & $0.013661344$ & $0.000119663$ & $-0.009406620$\\
8192 & $0.018392580$ & $0.000118390$ & $-0.004675384$\\\bottomrule
\end{tabular}
\caption{Computational bridge audit, not laboratory data. The million-path reference uses four independent seeds and paired nested meshes. The reported SE is the ratio influence-function estimate.}\label{tab:A}
\end{table}

Brownian-bridge conditioning to repair missed extrema is established methodology; see, for example, Shevchenko \cite{Shevchenko}. Here we repair the inferential target without pretending to have an exact joint range sampler. Conditional on neighboring bridge knots $u,v$ separated by $1/m$, a segment maximum has the law
\[
 \frac{u+v+\sqrt{(u-v)^2+2E/m}}2,\qquad E\sim\operatorname{Exp}(1).
\]
The corresponding minimum uses the negative square root. Independent arrays for maxima and minima give a synthetic range $\widetilde Y$ whose conditional mean equals that of $Y$, but whose conditional joint extrema are not exact. Appendix~\ref{app:A} proves a bias enclosure for smooth functions of this synthetic range and combines it with bounded empirical-Bernstein intervals. For the fixed design $m=1000$, $N=20000$, and upper clipping level $2.5$, the recorded computation gives
\[
 \widehat A(1)=0.02293831721,\qquad
 I_{0.95}=[0.01865431481,\ 0.02730880171].
\]
The coverage theorem is for ideal independent draws and exact arithmetic. The implementation uses pseudorandom numbers and floating-point operations, not a formally verified random-number program. This is a finite-sample, mesh-aware statistical construction, not evidence for the universal matrix hierarchy.

\section{Route B a finite-sample cumulant-gap observation}
Assume first that $N\geq2$ and that $X_1,\ldots,X_N$ are independent exact $q_0$ samples. The symmetry is known, so the estimator can use raw even moments without estimating a zero mean. Write
\[
 P(a,b,c)=ac+5a^2b-\tfrac{10}{3}b^2,
 \qquad (a,b,c)=(\mu_2,\mu_4,\mu_6).
\]
Its coordinates are strongly correlated. Treating three independently estimated cumulants as if their errors were independent is incorrect.

A transparent fixed-$N$ construction clips $|X|$ at a predetermined $L$ and builds simultaneous empirical-Bernstein intervals for the three bounded powers. An explicit kernel-tail bound supplies the one-sided clipping bias. The sharper version uses a fixed anchor $\theta_0=(a_0,b_0,c_0)$ and the scalar polynomial
\[
 h_0(x)=g_2x^2+g_4x^4+g_6x^6,\qquad
 g=(c_0+10a_0b_0,\ 5a_0^2-20b_0/3,\ a_0).
\]
For $d=(a,b,c)-\theta_0$, the exact Taylor remainder is
\begin{equation}\label{eq:remainder}
 P(\theta_0+d)-P(\theta_0)-g\cdot d
 =d_2d_6+10a_0d_2d_4+5b_0d_2^2+5d_2^2d_4-\tfrac{10}{3}d_4^2.
\end{equation}
Simultaneous moment intervals bound its absolute value by $R_{\max}$.
For the predetermined anchor $(0.04621,0.00596002,0.001205534)$ and $L=0.8$, $h_0$ is increasing in $|x|$ beyond $L$ and the clipped scalar has range less than $0.002944$. Thus clipping does not require a negative scalar-bias correction in a lower confidence bound. Define the scalar samples $H_i=h_0(\min(|X_i|,L))$. If $\bar h$ and $s_h^2$ are their sample mean and unbiased variance, put
\[
 e_h=\sqrt{\frac{2s_h^2\log(2/\alpha_h)}N}
      +\frac{7(0.002944)\log(2/\alpha_h)}{3(N-1)}.
\]
The empirical-Bernstein inequality \cite{MaurerPontil}, unioned with the moment events, gives
\begin{equation}\label{eq:Bbound}
 G_0\geq P(\theta_0)+\bar h-g\cdot\theta_0-e_h-R_{\max}
\end{equation}
with probability at least $1-\alpha_m-\alpha_h$. Appendix~\ref{app:B} gives the moment radii and bias constants. The same data may be used for all events; their independence is not required for the union bound. The anchor must be fixed before the test data, or chosen on an independent pilot split. This fixed-$N$ statement is not an optional-stopping guarantee.

A held-out numerical-law run with $N=10^6$, seed $2026110340$, and $\alpha_m=\alpha_h=0.025$ gave a raw gap $8.54804116\cdot10^{-7}$ and the ideal-exact-law lower endpoint $5.37384532\cdot10^{-7}$. The actual sampler was a floating-point inverse-CDF implementation. An interval audit supports a conservative CDF discrepancy envelope of $10^{-6}$ under its elementary-function, interpolation, floating-point, and random-input assumptions. Conditional on that discrepancy bound applying to the actual draw law, adding the deterministic bounded-variation transfer errors gives the corrected $95\%$ lower endpoint
\[
 5.19035463\cdot10^{-7}>0.
\]
This is not an unconditional bit-level certificate for the full generator. These bounds assume independent draws or independent preparations. A single-trap time series requires a dependence-aware confidence bound and an initialization analysis; stationarity or finite thinning alone does not ensure exact independence. More fundamentally, $G_0>0$ is already unconditional by Proposition~\ref{prop:signs}; measuring it validates the model or procedure, not RH.

\section{Discussion and limitations}
The model has a useful logical separation. The static susceptibility chain is already a classical GHS consequence, rather than a newly discovered implication of RH. The probability density, diffusion, susceptibility chain, finite-time response bound, deterministic-control bound, and quench theorem are unconditional. At $a=0$, all alternating zero-field signs and the particular low-order gap are also unconditional. Under RH one obtains a stronger positive Stieltjes representation and all finite matrix inequalities at once; at $a=0$ the universal L\"owner hierarchy is an equivalence. For $a\geq0.05$ that hierarchy is already unconditional, so observing it after enough quadratic deformation cannot adjudicate RH at the undeformed point.

The numerical work illustrates statistical and implementation constraints rather than replacing the analytic results. Bridge mesh bias, nonlinear ratio bias, correlated moment errors, numerical sampling-law discrepancy, finite feedback delay, and force saturation must each be accounted for in the observable to which they apply. A fluctuation measurement is not automatically a measurement of an infinite-order positivity property. A physical implementation would first require an independently calibrated force law, diffusivity, detector noise, time resolution, and model-discrepancy analysis.

The predicted variance quench is the central physical observable studied here, while the matrix hierarchy is an auxiliary mathematical diagnostic. The exact results suggest a tractable experimental program of model validation: test the force-dependent mean and variance, measure stationary step response, and compare a controlled quadratic quench with the monotone relaxation envelope. They do not suggest that finite-precision agreement can prove RH. Establishing new mathematical priority, obtaining a sharper implementation certificate, and demonstrating the potential in a laboratory remain outside the claims of this paper.

\paragraph{Data and code availability.}
The accompanying reproducibility package contains the scripts and machine-readable outputs underlying the reported computations, with an inventory distinguishing exact analytic arguments, interval enclosures, ordinary quadrature/PDE calculations, and pseudorandom simulations. This manuscript does not claim a new physical experiment, public dataset deposit, journal submission, or external peer review. Randomized coverage statements specify their ideal independent-input assumptions. The numerical code is not proof-assistant verified.

\clearpage
\section*{AI statement}
Under the author's direction, AI (ChatGPT, OpenAI; version not recorded) assisted with mathematical exploration, numerical scripts, drafting, and editing. The preparation checks are described in the methods. Independent verification by the author remains necessary before submission; the author retains responsibility for the final content.

\section*{Appendix}
\appendix
\section{Full proof of the delayed-feedback bound}\label{app:delayproof}

Let $b=b_2$ be the smooth ideal drift and let
\begin{equation}\label{eq:ideal}
 \dd X_t=b(X_t)\dd t+\dd W_t,\qquad X_0\sim q_0.
\end{equation}
The initial variable is independent of the standard Brownian motion $W$.
Write $t_n=n\Delta$ and $X_n=X_{t_n}$. Set
\[
 g(x)=\clip(b(x),-B,B),\qquad
 G(x)=\E_\xi[g(x+\sigma\xi)],\qquad
 K=\Lip(g),\qquad \xi\sim N(0,1).
\]
Here $g$ is nonincreasing and globally $K$-Lipschitz, so $G$ has both
properties as well. Let $(\xi_j)_{j\geq-1}$ be independent standard normals,
independent jointly of $X_0$ and the entire Brownian motion. The implemented
process, sampled at the update times, obeys exactly
\begin{equation}\label{eq:scheme}
 Z_{n+1}=Z_n+\Delta g(Z_{n-1}+\sigma\xi_{n-1})
                 +(W_{t_{n+1}}-W_{t_n}),\qquad n\geq0,
\end{equation}
where $Z_{-1}=Z_0=X_0$ and the formal ideal value is $X_{-1}=X_0$.
The same stored $\xi_{n-1}$ is used when its measurement is actuated;
\eqref{eq:scheme} is not a fresh noisy observation of $Z_n$.

The ideal law is dominated at every nonnegative time by
\begin{equation}\label{eq:domination}
 \mathcal L(X_t)(\dd x)\leq\rho(x)\dd x,
 \qquad \rho(x)=q_0(x)e^{2x^2}.
\end{equation}
For completeness, the invariant density is $q_2=\rho/Z(2,0)$.
Initially $q_0\leq\rho$, and positivity of the Markov transition kernel,
together with invariance of the finite measure $\rho(x)\dd x$, gives
\eqref{eq:domination}. Define
\begin{equation}\label{eq:constants}
 \begin{split}
 a(x)&=(Lb)(x)=b(x)b'(x)+\tfrac12 b''(x),\\
 C_1^2&=\int_{\R}\rho(x)b'(x)^2\dd x,\qquad
 C_0^2=\int_{\R}\rho(x)a(x)^2\dd x,\\
 C_s^2&=\int_{\R}\rho(x)(G(x)-b(x))^2\dd x.
 \end{split}
\end{equation}
These constants are finite in the theta-kernel model. In particular,
\begin{equation}\label{eq:momentbounds}
 \norm{b'(X_t)}_2\leq C_1,\quad
 \norm{a(X_t)}_2\leq C_0,\quad
 \norm{G(X_t)-b(X_t)}_2\leq C_s
 \qquad (t\geq0).
\end{equation}
The regularity and integrability in these statements justify It\^o's formula
and the ordinary and stochastic Fubini steps used below; equivalently one
can first localize the smooth diffusion and pass to the limit using these
square-integrability bounds.

\begin{proposition}\label{prop:expanded}
Let $N\geq1$ and $\alpha>0$. If
\begin{equation}\label{eq:stepsize}
 \Delta K\leq \frac{2\alpha}{\alpha^2+\alpha+1},
\end{equation}
then
\begin{equation}\label{eq:bound}
 \boxed{
 \norm{Z_N-X_N}_2\leq\sqrt{1+\alpha}
 \left\{
 \left[\frac{7N-6}{3}C_1^2\Delta^3
             +N\Delta^2K^2\sigma^2\right]^{1/2}
 +\frac{3N-2}{2}\Delta^2 C_0+N\Delta C_s
 \right\}.}
\end{equation}
This is a bound for the point-observation, one-step-delayed model with the
specified exact common initialization. It does not require RH.
\end{proposition}

\begin{proof}
\par\noindent\textbf{1. Two filtrations, with distinct purposes.} 
For $n\geq0$ put
\begin{equation}\label{eq:gridfiltration}
 \mathcal F_n=
 \sigma\left(X_0,\,W_s:0\leq s\leq t_n,\,
                  \xi_j:-1\leq j\leq n-2\right),
\end{equation}
completed by null sets. The noise set in $\mathcal F_0$ is empty.
Both $X_n$ and $Z_n$ are $\mathcal F_n$-measurable: at the end of step
$n-1$, the most recently actuated measurement was $\xi_{n-2}$.
Thus $\xi_{n-1}$ is independent of $\mathcal F_n$.
This is a reduced filtration for the proof. A measurement may physically
be stored earlier, but the state does not use it before its actuation, so
we may leave it unrevealed in \eqref{eq:gridfiltration} until that time.

We also use the usual augmentation of
\begin{equation}\label{eq:brownianfiltration}
 \mathcal H_t=
 \sigma\left(X_0,\,W_s:0\leq s\leq t,\,
                         \xi_j:j\geq-1\right).
\end{equation}
All measurement noises are revealed at time zero in this second
filtration. Their joint independence from $W$ implies that $W$ remains
Brownian with respect to $\mathcal H_t$. In particular, square-integrable
It\^o integrals over $[u,v]$ have conditional mean zero given
$\mathcal H_u$. We use $\mathcal F_n$ to center a fresh measurement error,
and $\mathcal H_t$ to center a Brownian integral even when an independent
measurement noise has already been included in its multiplier.

\par\noindent\textbf{2. Exact local It\^o residual.} 
For $n\geq0$, define
\begin{equation}\label{eq:MQ}
 \begin{split}
 Q_n&=\int_{t_n}^{t_{n+1}}b'(X_s)\dd W_s,\\
 M_n&=\int_{t_n}^{t_{n+1}}(t_{n+1}-s)b'(X_s)\dd W_s,
 \end{split}
 \qquad Q_{-1}=0.
\end{equation}
It\^o's formula for the smooth function $b$ gives
\[
 b(X_s)-b(X_u)=\int_u^s a(X_v)\dd v
                         +\int_u^s b'(X_v)\dd W_v.
\]
At the first step, ordinary and stochastic Fubini therefore give
\begin{equation}\label{eq:residual0}
 \begin{split}
 X_1-X_0-\Delta b(X_{-1})-(W_{t_1}-W_{t_0})
 &=\int_{t_0}^{t_1}[b(X_s)-b(X_0)]\dd s\\
 &=R_0+M_0,\\
 R_0&=\int_{t_0}^{t_1}(t_1-s)a(X_s)\dd s.
 \end{split}
\end{equation}
For $n\geq1$, separate the preceding complete interval from the current
partial interval:
\begin{align*}
 &\int_{t_n}^{t_{n+1}}[b(X_s)-b(X_{n-1})]\dd s\\
 &\quad=\Delta[b(X_n)-b(X_{n-1})]
          +\int_{t_n}^{t_{n+1}}[b(X_s)-b(X_n)]\dd s\\
 &\quad=\Delta\int_{t_{n-1}}^{t_n}a(X_s)\dd s
         +\Delta Q_{n-1}
         +\int_{t_n}^{t_{n+1}}(t_{n+1}-s)a(X_s)\dd s+M_n.
\end{align*}
Consequently the exact formula for every $n\geq0$ is
\begin{equation}\label{eq:localresidual}
 X_{n+1}-X_n-\Delta b(X_{n-1})-(W_{t_{n+1}}-W_{t_n})
       =R_n+M_n+\Delta Q_{n-1},
\end{equation}
where, for $n\geq1$,
\begin{equation}\label{eq:residualn}
 R_n=\Delta\int_{t_{n-1}}^{t_n}a(X_s)\dd s
                +\int_{t_n}^{t_{n+1}}(t_{n+1}-s)a(X_s)\dd s.
\end{equation}
Minkowski's integral inequality and \eqref{eq:momentbounds} give
\begin{equation}\label{eq:Rbounds}
 \norm{R_0}_2\leq\frac{\Delta^2}{2}C_0,
 \qquad
 \norm{R_n}_2\leq\frac{3\Delta^2}{2}C_0\quad(n\geq1).
\end{equation}
The stochastic Fubini integrals are square integrable. For example,
\[\int_{t_n}^{t_{n+1}}(t_{n+1}-s)^2\E b'(X_s)^2\dd s
\leq C_1^2\Delta^3/3.\]
No centeredness is asserted for $R_n$.

\par\noindent\textbf{3. Exact error recursion and the measurement residual.} 
Write $E_n=Z_n-X_n$, so $E_{-1}=E_0=0$. For $n\geq1$, choose
\begin{equation}\label{eq:c}
 c_n=
 \begin{cases}
 -\Delta\dfrac{G(Z_{n-1})-G(X_{n-1})}{Z_{n-1}-X_{n-1}},
                                     & E_{n-1}\ne0,\\[2mm]
 0,                                  & E_{n-1}=0,
 \end{cases}
 \qquad c_0=0.
\end{equation}
Then $0\leq c_n\leq\Delta K$ and $c_n$ is
$\mathcal F_{n-1}$-measurable for $n\geq1$. In particular it is determined
\emph{before} the Brownian interval $[t_{n-1},t_n]$.
The zero-error convention leaves the exact identity
$\Delta[G(Z_{n-1})-G(X_{n-1})]=-c_nE_{n-1}$ unchanged.

Define
\begin{equation}\label{eq:eta}
 \eta_n=\Delta\{g(Z_{n-1}+\sigma\xi_{n-1})-G(Z_{n-1})\}.
\end{equation}
Since $Z_{n-1}$ is $\mathcal F_n$-measurable, the fresh-noise independence
in Step~1 proves
\begin{equation}\label{eq:etamean}
 \E[\eta_n\mid\mathcal F_n]=0.
\end{equation}
For a fixed $z$, and independent standard normals $\xi,\xi'$, the elementary
independent-copy variance identity gives
\begin{align*}
 \Var(g(z+\sigma\xi))
 &=\tfrac12\E[(g(z+\sigma\xi)-g(z+\sigma\xi'))^2]\\
 &\leq\tfrac12K^2\sigma^2\E[(\xi-\xi')^2]
 =K^2\sigma^2.
\end{align*}
Conditioning on $\mathcal F_n$ therefore yields
\begin{equation}\label{eq:etavariance}
 \E[\eta_n^2\mid\mathcal F_n]\leq\Delta^2K^2\sigma^2.
\end{equation}
Subtracting \eqref{eq:localresidual} from \eqref{eq:scheme} now gives
\begin{equation}\label{eq:errorrecursion}
 E_{n+1}=E_n-c_nE_{n-1}+U_n+r_n,
\end{equation}
where
\begin{equation}\label{eq:Ur}
 U_n=-M_n-\Delta Q_{n-1}+\eta_n,
 \qquad
 r_n=\Delta[G(X_{n-1})-b(X_{n-1})]-R_n.
\end{equation}
The formal $X_{-1}=X_0$ is used in $r_0$. Equations
\eqref{eq:momentbounds} and \eqref{eq:Rbounds} imply
\begin{equation}\label{eq:rsum}
 \sum_{n=0}^{N-1}\norm{r_n}_2
 \leq N\Delta C_s+
       \left(\tfrac12+\tfrac32(N-1)\right)\Delta^2C_0
 =N\Delta C_s+\frac{3N-2}{2}\Delta^2C_0.
\end{equation}
Crucially, $U_n$ is \emph{not} asserted to be a martingale difference:
$Q_{n-1}$ is already $\mathcal F_n$-measurable.

\par\noindent\textbf{4. A common pathwise contracting norm.} 
For $u,v\in\R$ define
\[
 \norm{(u,v)}_P^2=u^2+\alpha(u-v)^2,
 \qquad
 P=\begin{pmatrix}1+\alpha&-\alpha\\-\alpha&\alpha\end{pmatrix},
 \qquad
 A(c)=\begin{pmatrix}1&-c\\1&0\end{pmatrix}.
\]
The matrix $P$ is positive definite since $\alpha>0$.
In the coordinates $(u,d)$, where $d=u-v$, $A(c)$ sends
\[
 (u,d)\longmapsto ((1-c)u+cd,-cu+cd).
\]
Thus the difference between the old and new squared norms is the quadratic
form with matrix
\begin{equation}\label{eq:decrement}
 D(c)=
 \begin{pmatrix}
 2c-(1+\alpha)c^2 & -c+(1+\alpha)c^2\\
 -c+(1+\alpha)c^2 & \alpha-(1+\alpha)c^2
 \end{pmatrix}.
\end{equation}
Its determinant is exactly
\begin{equation}\label{eq:det}
 \det D(c)=c\{2\alpha-(\alpha^2+\alpha+1)c\}.
\end{equation}
If $0<c\leq2\alpha/(\alpha^2+\alpha+1)$, the first diagonal entry is
strictly positive, because
\[
 \frac{2\alpha}{\alpha^2+\alpha+1}<\frac{2}{1+\alpha}.
\]
The determinant in \eqref{eq:det} is nonnegative, so the second diagonal
entry is nonnegative as well and $D(c)$ is positive semidefinite. At $c=0$,
$D(0)=\operatorname{diag}(0,\alpha)$ is also positive semidefinite.
Consequently \eqref{eq:stepsize} gives, for every realization and every $n$,
\begin{equation}\label{eq:contract}
 \norm{A(c_n)(u,v)}_P\leq\norm{(u,v)}_P.
\end{equation}
An input in the first coordinate has norm
$\norm{(w,0)}_P=\sqrt{1+\alpha}\,|w|$, and
$|u|\leq\norm{(u,v)}_P$.

\par\noindent\textbf{5. Algebraic split with the realized coefficients.} 
Define $S_{-1}=S_0=D_{-1}=D_0=0$ and
\begin{equation}\label{eq:split}
 \begin{split}
 S_{n+1}&=S_n-c_nS_{n-1}+U_n,\\
 D_{n+1}&=D_n-c_nD_{n-1}+r_n.
 \end{split}
\end{equation}
The coefficients in these two equations are exactly the realized $c_n$
from \eqref{eq:c}; they are not recomputed from $S$ or from $D$.
Induction in \eqref{eq:errorrecursion} proves $E_n=S_n+D_n$.
The split is an algebraic identity even though the coefficients depend on
the full error process. We neither condition on the entire coefficient
sequence nor claim that it is independent of the forcing.
Induction also shows that $S_n$ and $D_n$ are $\mathcal F_n$-measurable.

By pathwise contraction and the triangle inequality,
\[
 \norm{(D_N,D_{N-1})}_P
 \leq\sqrt{1+\alpha}\sum_{n=0}^{N-1}|r_n|.
\]
Applying the $L^2$ triangle inequality and \eqref{eq:rsum} gives
\begin{equation}\label{eq:Denergy}
 \left(\E\norm{(D_N,D_{N-1})}_P^2\right)^{1/2}
 \leq\sqrt{1+\alpha}
       \left(N\Delta C_s+\frac{3N-2}{2}\Delta^2C_0\right).
\end{equation}
This argument does not require $r_n$ to be adapted before its own interval,
centered, or independent of any coefficient.

\par\noindent\textbf{6. Every stochastic energy cross term.} 
Put
\[
 V_n=\norm{(S_n,S_{n-1})}_P^2,\qquad
 \ell_n=(1+\alpha)(S_n-c_nS_{n-1})-\alpha S_n
        =S_n-(1+\alpha)c_nS_{n-1}.
\]
The second coordinate of $A(c_n)(S_n,S_{n-1})$ is $S_n$, not
$S_{n-1}$; this determines the $-\alpha S_n$ term in $\ell_n$.
Equivalently, direct expansion of
$(S_n-c_nS_{n-1}+U_n)^2+\alpha(-c_nS_{n-1}+U_n)^2$
gives the same multiplier. Expanding the squared norm in the first equation
of \eqref{eq:split} gives the exact identity, followed by
\eqref{eq:contract},
\begin{align}
 V_{n+1}
 &=\norm{A(c_n)(S_n,S_{n-1})}_P^2
                    +2\ell_nU_n+(1+\alpha)U_n^2\nonumber\\
 &\leq V_n+2\ell_nU_n+(1+\alpha)U_n^2.
 \label{eq:energyidentity}
\end{align}
The multiplier $\ell_n$ is $\mathcal F_n$-measurable.
For the current measurement and current Brownian integral this implies
\begin{equation}\label{eq:zerocurrent}
 \E[\ell_n\eta_n]=0,\qquad \E[\ell_nM_n]=0.
\end{equation}
The first equality follows by conditioning on $\mathcal F_n$ and using
\eqref{eq:etamean}; the second follows by conditioning on
$\mathcal H_{t_n}$ and using the Brownian martingale property.

The preceding Brownian interval must be treated differently. For $n\geq1$,
substitution of the recurrence for $S_n$ gives
\begin{equation}\label{eq:ell}
 \ell_n=B_n-M_{n-1},
\end{equation}
where
\begin{equation}\label{eq:B}
 \begin{split}
 B_n={}&[1-(1+\alpha)c_n]S_{n-1}
       -c_{n-1}S_{n-2}\\
       &-\Delta Q_{n-2}+\eta_{n-1}.
 \end{split}
\end{equation}
For $n=1$ this uses $c_0=0$, $S_{-1}=0$, and $Q_{-1}=0$.
Each term of $B_n$ is $\mathcal H_{t_{n-1}}$-measurable:
$c_n$ and $S_{n-1}$ are known by $t_{n-1}$; $Q_{n-2}$ is an integral
ending then; and $\eta_{n-1}$ uses only $Z_{n-2}$ and $\xi_{n-2}$,
both included in $\mathcal H_{t_{n-1}}$.
Therefore
\[
 \E[B_nQ_{n-1}]=
 \E\{B_n\E[Q_{n-1}\mid\mathcal H_{t_{n-1}}]\}=0.
\]
It follows that the one surviving cross term is exactly
\begin{equation}\label{eq:survivor}
 \E[\ell_n(-\Delta Q_{n-1})]
  =\Delta\E[M_{n-1}Q_{n-1}].
\end{equation}
At $n=0$ there is no preceding interval and $\ell_0=Q_{-1}=0$.
In particular, discarding \eqref{eq:survivor} would be incorrect.

For completeness, the three distinct components of $U_n$ are pairwise
orthogonal in expectation. Indeed,
\begin{equation}\label{eq:Uorthogonal}
 \begin{aligned}
 \E[M_nQ_{n-1}]&=0
   &&\text{by conditioning on }\mathcal H_{t_n},\\
 \E[\eta_nQ_{n-1}]&=0
   &&\text{by conditioning on }\mathcal F_n,\\
 \E[\eta_nM_n]&=0
   &&\text{by conditioning on }\mathcal H_{t_n}.
 \end{aligned}
\end{equation}
For the last equality, $\eta_n$ is
$\mathcal H_{t_n}$-measurable, even though it is not
$\mathcal F_n$-measurable. Hence
\begin{equation}\label{eq:Uvariance}
 \E U_n^2=\E M_n^2+\Delta^2\E Q_{n-1}^2+\E\eta_n^2.
\end{equation}
Equations \eqref{eq:zerocurrent}, \eqref{eq:survivor}, and
\eqref{eq:Uvariance} account for every cross term in
\eqref{eq:energyidentity}.

\par\noindent\textbf{7. Summation and the constants $7/3$ and $1/3$.} 
Taking expectations in \eqref{eq:energyidentity}, summing from $n=0$ to
$N-1$, and using $V_0=0$ gives
\begin{align}
 \E V_N\leq{}&(1+\alpha)\left\{
 \sum_{n=0}^{N-1}\E M_n^2
 +\Delta^2\sum_{n=1}^{N-1}\E Q_{n-1}^2
 +\sum_{n=0}^{N-1}\E\eta_n^2\right\}\nonumber\\
 &\quad+2\Delta\sum_{n=1}^{N-1}\E[M_{n-1}Q_{n-1}].
 \label{eq:energysum}
\end{align}
The exact covariance coefficient in this formula is $2\Delta$, not
$2(1+\alpha)\Delta$. It\^o covariance isometry gives
\begin{equation}\label{eq:covpositive}
 \E[M_jQ_j]=\int_{t_j}^{t_{j+1}}(t_{j+1}-s)
                                    \E b'(X_s)^2\dd s\geq0.
\end{equation}
We may consequently \emph{increase} the right-hand side of
\eqref{eq:energysum} by replacing $2\Delta$ with
$2(1+\alpha)\Delta$. Reindexing $j=n-1$ and grouping the resulting
squares gives the valid, slightly conservative estimate
\begin{equation}\label{eq:energygrouped}
 \E V_N\leq(1+\alpha)\left\{
 \sum_{j=0}^{N-2}\E(M_j+\Delta Q_j)^2
 +\E M_{N-1}^2+\sum_{n=0}^{N-1}\E\eta_n^2\right\}.
\end{equation}
When $N=1$, the first sum is empty. The grouping is performed only
\emph{after} taking the adapted energy expectations, so it does not
replace the delayed forcing by an incorrectly independent sequence.

By \eqref{eq:MQ}, the combined integral is
\[
 M_j+\Delta Q_j
 =\int_{t_j}^{t_{j+1}}(t_{j+1}-s+\Delta)b'(X_s)\dd W_s.
\]
It\^o isometry and \eqref{eq:momentbounds} therefore give
\begin{align}
 \E(M_j+\Delta Q_j)^2
 &=\int_{t_j}^{t_{j+1}}(t_{j+1}-s+\Delta)^2
                                \E b'(X_s)^2\dd s\nonumber\\
 &\leq C_1^2\int_0^\Delta(2\Delta-u)^2\dd u
 =\frac73 C_1^2\Delta^3,
 \label{eq:seven}\\
 \E M_{N-1}^2
 &\leq C_1^2\int_0^\Delta(\Delta-u)^2\dd u
 =\frac13 C_1^2\Delta^3.
 \label{eq:one}
\end{align}
The covariance in \eqref{eq:covpositive} has a positive sign in the
energy bound, not a cancellation to be set to zero. The improvement is that all these local contributions accumulate
in squared energy, rather than by summing their individual $L^2$ norms.
Using \eqref{eq:etavariance} in \eqref{eq:energygrouped} now gives
\begin{equation}\label{eq:Senergy}
 \E\norm{(S_N,S_{N-1})}_P^2
 \leq(1+\alpha)\left\{
 \left(\frac73(N-1)+\frac13\right)C_1^2\Delta^3
                      +N\Delta^2K^2\sigma^2\right\}.
\end{equation}
The coefficient in parentheses is $(7N-6)/3$.

\par\noindent\textbf{8. Final $L^2$ bound.} 
Since $E_N=S_N+D_N$, the triangle inequality in
$L^2(\Omega;\R^2,P)$ gives
\begin{align*}
 \norm{E_N}_2
 &\leq\left(\E\norm{(E_N,E_{N-1})}_P^2\right)^{1/2}\\
 &\leq\left(\E\norm{(S_N,S_{N-1})}_P^2\right)^{1/2}
       +\left(\E\norm{(D_N,D_{N-1})}_P^2\right)^{1/2}.
\end{align*}
Substitution of \eqref{eq:Denergy} and \eqref{eq:Senergy} proves
\eqref{eq:bound}.
\end{proof}

\section{Integral tails and interval reproduction}\label{app:integrals}
Let $I_k=\int_0^\infty x^kK(x)\dd x$, so $\mu_k=I_k/I_0$ for even $k$. The interval scripts use the first three terms on $[0,2]$ and explicit tails. Here is an elementary majorant sufficient for the reported design. For the series tail $n\geq4$, $x\in[0,2]$, use $\pi>3$, $\pi<4$, $e<3$, and $2R_n-3<2R_n$ to obtain
\[
 \sum_{n\geq4}e^{x/2}R_n(2R_n-3)e^{-R_n}
 \leq B e^{-91x},\qquad
 B=\frac{8192e^{-48}}{1-e^{-23}}.
\]
For clarity, $R_n\geq3n^2(1+2x)$ and the polynomial prefactor is bounded by $32n^4e^{(9/2)x}$; at $n=4$ this gives $8192e^{-48}e^{-(96-9/2)x}$, which is at most $8192e^{-48}e^{-91x}$. Successive $n$ terms have ratio at most $(5/4)^4e^{-27}<e^{-23}$, uniformly in $x\geq0$. Thus the $k$th-moment series remainder is at most $Bk!/91^{k+1}$.

For the spatial tail $x\geq2$, write $u=x-2$ and use $\pi e^4>160$. One may bound the full kernel by
\[
 K(2+u)\leq C e^{-315u},\qquad
 C=\frac{6\cdot160^2e^{-160}}{1-e^{-316}}.
\]
To verify the prefactor, at $x=2$ one has $e^{x/2}=e<3$, and $r^2e^{-r}$ decreases for $r\geq2$; the first term is therefore at most $6\cdot160^2e^{-160}$. For each term, the logarithmic derivative is $9/2+6/(2R_n-3)-2R_n\leq9/2+6/317-320<-315$, and the remaining $n$ terms are bounded geometrically. Hence
\[
 \int_2^\infty x^kK(x)\dd x
 \leq C\sum_{j=0}^k\binom{k}{j}2^{k-j}\frac{j!}{315^{j+1}}.
\]
For the cosine integral the absolute tail bound is simply $B/91+C/315$.
These bounds intentionally sacrifice sharpness. Directed interval evaluation of the boxes and all constants is what makes the enclosure, rather than apparent decimal stability, the relevant numerical object.

\section{A mesh-aware finite-sample bridge construction}\label{app:A}
Conditional on the $m+1$ bridge knots, residual segment bridges are independent. For the exact range and the synthetic range described in Section~6,
\[
 \E(Y\mid\mathrm{knots})=\E(\widetilde Y\mid\mathrm{knots}),\qquad
 \Var(Y\mid\mathrm{knots}),\ \Var(\widetilde Y\mid\mathrm{knots})
 \leq\frac1{\pi m}.
\]
To see the variance bound, the conditional residual at any point has variance at most $1/(4m)$, and conditional covariances within a segment are nonnegative. The squared Gaussian gradient norm of the range is at most $1/(2m)$; independent upper and lower residual copies obey the same bound. Gaussian Poincar\'e and the factor $2/\pi$ prove the assertion. Suprema can be approximated on finite dense grids and passed to the limit in $L^2$.

If $H$ is convex with $0\leq H''\leq K_H$, its conditional Jensen gap lies in $[0,K_H/(2\pi m)]$ for either variable. Their means agree, so
\begin{equation}\label{eq:Abias}
 |\E H(\widetilde Y)-\E H(Y)|\leq\frac{K_H}{2\pi m}.
\end{equation}
There is no extra factor of two, since two numbers in the same interval differ by at most the interval length.

Take $f(y)=y^{3/2}\log y$ and $g(y)=y^{3/2}$. Extend each affinely below $a_*=1/4$ by its tangent there, obtaining $f_*,g_*$. Their second derivatives are nonnegative and bounded by $K_f=1.5e^{1/3}$ and $K_g=1.5$. Their absolute extension errors below $a_*$ are less than one. The exact bridge-range distribution \cite{PermanWellner} gives
\[
 P(Y\leq a_*)=4\pi a_*^{-3}\sum_{k\geq1}k^2e^{-\pi k^2/a_*^2}
 <1.2\cdot10^{-19}.
\]
This controls the extension error for the target; no claim about the small tail of $\widetilde Y$ is needed because the estimator actually uses the extensions.

For the upper tail, a conditional independent copy of the negative minimum gives a coupling
$\widetilde Y=Y+\sqrt{2/\pi}(A'-A)$. Gaussian concentration yields, for $0<\alpha<1$,
\[
 P(\widetilde Y>x)\leq P(Y>\alpha x)
       +\exp\{-\pi m(1-\alpha)^2x^2/2\}.
\]
For $\alpha=0.95$ and $x\geq L=2.5$, the exact survival series gives
$P(Y>\alpha x)\leq13\alpha^2x^2e^{-\pi\alpha^2x^2}$, and $f'(x),g'(x)\leq x^2$. Therefore clipping at $L$ costs at most
\[
 T_m=13\alpha^2 I_4(\pi\alpha^2)
        +I_2(\pi m(1-\alpha)^2/2),\qquad
 I_j(c)=\int_L^\infty x^je^{-cx^2}\dd x.
\]
The recurrence
$I_j=L^{j-1}e^{-cL^2}/(2c)+(j-1)I_{j-2}/(2c)$ and
$I_0\leq e^{-cL^2}/(2cL)$ give elementary numerical upper bounds. At $m=1000$, $T_m<7.077770\cdot10^{-7}$. Thus total expectation biases are bounded by
$b_f<0.000333886$ and $b_g<0.000239441$.

For independent synthetic samples use $U_i=f_*(\min(\widetilde Y_i,L))$ and $D_i=g_*(\min(\widetilde Y_i,L))$. Their range lengths are at most
$R_U=L^{3/2}\log L+2/(3e)$ and $R_D=L^{3/2}+1/16$. With unbiased sample variances $s_H^2$, set
\[
 r_H=\sqrt{2s_H^2\log(8/\delta)/N}
       +\frac{7R_H\log(8/\delta)}{3(N-1)},\qquad H=U,D.
\]
Maurer--Pontil's bounded empirical-Bernstein inequality and a union bound give simultaneous intervals
$[\bar H-r_H-b_H,\bar H+r_H+b_H]$ for the target numerator and denominator, with probability at least $1-\delta$. When the denominator interval is positive, the minimum and maximum of the four corner ratios $U/(2D)$ enclose $A(1)$. This proves the stated fixed-design confidence guarantee.

\section{Elementary bounds for the fast-feedback design}\label{app:fast}
For $x\geq0$ write $R=\pi e^{2x}$, $R_n=n^2R$, and $f_n=e^{x/2}R_n(2R_n-3)e^{-R_n}$, with weights $w_n=f_n/\sum_j f_j$. Their logarithmic derivatives are
\[
 p_n=\frac92+\frac6{2R_n-3}-2R_n,\qquad
 p_n'=-4R_n-\frac{24R_n}{(2R_n-3)^2}.
\]
For $b=b_0+2x$, $b'=2+(\sum_n w_np_n'+\operatorname{Var}_w p_n)/2$. Using the already proved $b'<-7$, positivity of the variance and decreasing $r/(2r-3)^2$ for $r\geq\pi$ give
\[
 |b'|\leq 2\E_w R_n+12\E_w\frac{R_n}{(2R_n-3)^2}-2.
\]
The ratio estimate $f_n/f_1\leq [2\pi/(2\pi-3)]n^4e^{-(n^2-1)R}$ implies $\E_w R_n/R<1.01$, by isolating $n=2$ and bounding the remaining series geometrically. The same positive-series bounds give $b(1.2)>-30>b(1.25)$. Thus the clip breakpoint lies in $(1.2,1.25)$ and
\[
 K_{30}\leq 2(1.01)\pi e^{2.5}+\frac{12\pi}{(2\pi-3)^2}-2<80.
\]
The clipped drift is constant outside the breakpoint; it is not strongly dissipative there.

For the tail, $R_0=\pi e^{2.4}$ and the series bounds give $|b(x)|\leq1.05R$, $K(x)\leq2.004\pi^{-1/4}R^{9/4}e^{-R}$ for $x\geq1.2$. The normalization exceeds $0.12$: keep the first positive term on $[-0.25,0.25]$. Put
\[
 \eta=\frac{\tfrac12\log^2(R_0/\pi)}{R_0},\quad c_*=1-\eta,\quad
 C=\frac{2.004\pi^{-1/4}}{0.12},\quad
 A_0=1.05R_0-30,\quad \lambda=c_*-\frac5{4R_0}>0.
\]
Since $\log^2(R/\pi)/R$ decreases for $R\geq R_0$, the inequality $e^{2x^2}\leq e^{\eta R}$ holds. Substitution $R=\pi e^{2x}$ and symmetry, followed by $(R_0+u)^{5/4}\leq R_0^{5/4}e^{5u/(4R_0)}$, yield the elementary expression
\[
 \epsilon_{30}^2\leq C R_0^{5/4}e^{-c_*R_0}
 \left(\frac{A_0^2}{\lambda}+\frac{2(1.05)A_0}{\lambda^2}+\frac{2(1.05)^2}{\lambda^3}\right)
 <(3.4\cdot10^{-5})^2.
\]
The last decimal inequality uses ordinary arithmetic with margin. The accompanying script evaluates this explicit expression; it does not certify floating-point rounding.

\section{Constants in the one-step-delay bound}\label{app:delayconstants}
The code encloses the first three positive theta terms and their logarithmic derivatives through third order on $15800$ cells of width $10^{-4}$. It bounds the omitted $n\geq4$ relative derivative tails by $10^{-10}$ on the whole range, and supplies explicit spatial tails beyond $|x|=1.5$. Upper and lower cell sums enclose the normalization and the moments; interval bounds on $b'$, $bb'+b''/2$, and shifted $b''$ then give the displayed $C_1,C_0$ bounds. The computation and its analytic remainders are detailed in the accompanying \texttt{interval\_constants\_proof.md} and script. They rely on the interval library's elementary enclosures and the explicitly inflated floating-point aggregation, not on bit-level proof-assistant verification.

For the Lipschitz bound, the positive series gives $b(1.21)<-30.615$ and $b(1.2)>-30$, so the clipping breakpoint is below $1.21$. The same series estimate as in Appendix~\ref{app:passive} bounds $\E_w R_n-R<0.001$ for $x\geq1$. Newman's convexity of $V'$ makes $-b'$ nondecreasing for $x\geq0$. With $R_*=\pi e^{2.42}$,
\[
 K\leq2R_*+\frac{12R_*}{(2R_*-3)^2}+0.002-2<70.
\]
For $C_s$, split the Gaussian noise into $|\xi|\leq8$ and its tail. A symmetric Taylor expansion bounds smoothing by $\sigma^2C_2/2$, with $C_2$ an $L^2(\rho)$ supremum of $|b''(x+u)|$ over $|u|\leq0.08$. Add the analogous shifted clipping-tail norm and a bound below $10^{-12}$ for $|\xi|>8$. The shifted tail bounds are checked down to $|x|=1.42$ and the script verifies $\int\rho<2$, which also bounds the Gaussian-tail contribution. The resulting conservative upper bound is $C_s<0.00135$.

\section{Bounds for the passive-residual design}\label{app:passive}
For $x\geq1$, put $\rho=\pi e^{2x}$ and $\rho_1=\pi e^2$. As in Appendix~\ref{app:fast}, the positive summands satisfy $f_n/f_1\leq2n^4e^{-(n^2-1)\rho}$. Their weighted correction to $\rho$ is below $0.001$. Newman proved convexity of $V'$, where $V=-\log q_0$ \cite{Newman1991}, so $-r'(x)$ is nondecreasing on $[0,1]$. The logarithmic derivative formulas in Appendix~\ref{app:fast} then give
\[
 |r(1)|\leq \rho_1-11.25+0.001<12,\quad
 \operatorname{Lip}(r_1)\leq 2\rho_1+\frac{12\rho_1}{(2\rho_1-3)^2}+0.002-9<38.
\]
Here the logarithmic-derivative variance term only decreases the upper bound on $-b_0'$. For $x\geq1$, the same series bound gives $-r'(x)\leq2\rho$ and hence $|r(x)-r(1)|\leq\rho-\rho_1$. In addition $K(x)\leq3e^{x/2}\rho^2e^{-\rho}$. The interval-box computation establishes total kernel mass $>0.24$. Writing $\epsilon_{\rm tail}^2:=\int q_0(x)e^{2x^2}|r(x)-r_1(x)|^2\dd x$, we obtain
\[
 \epsilon_{\rm tail}^2\leq25\int_1^\infty e^{x/2+2x^2}\rho^2(\rho-\rho_1)^2e^{-\rho}\dd x.
\]
Set $x=1+y$; use $e^{2y}-1\leq2ye^{2y}$ and $e^{2y}\geq1+2y+2y^2$ to obtain
\[
 \epsilon_{\rm tail}^2\leq\frac{200\rho_1^4e^{2.5-\rho_1}}{(2\rho_1-12.5)^3}
 \leq\frac{200\cdot23^4e^{2.5-23}}{(46-12.5)^3}<0.002^2.
\]
Thus we may take $\epsilon=0.002$ in Proposition~\ref{prop:passive}. The last upper bound uses monotonicity in $\rho_1\geq23$. The reproducibility package includes the interval boxes, tail majorants, and numerical endpoint evaluation; its enclosures rely on the stated interval-library assumptions.

\section{Moment clipping and the covariance-aware gap bound}\label{app:B}
The density induced by \eqref{eq:sizebias}, for $x\geq0$, may also be written as
\[
 q_0(x)=Z_*^{-1}\sum_{n\geq1}
 [8\pi^2n^4e^{9x/2}-12\pi n^2e^{5x/2}]e^{-\pi n^2e^{2x}},
 \qquad Z_*=\E_PY^{1/2}.
\]
The range dominates the bridge maximum and $P(\max B>r)=e^{-2r^2}$; hence
$Z_*\geq\Gamma(5/4)/\pi^{1/4}=:\underline Z$.
For a clipping level $L\geq0$, put
\begin{align*}
 A_L=\pi e^{2L},\quad r_L=e^{-A_L},\quad
 C_L&=\frac{1+11r_L+11r_L^2+r_L^3}{(1-r_L)^5},\\
 D_L=\frac{8\pi^2}{\underline Z}C_Le^{9L/2-A_L},\quad
 \lambda_L&=2A_L-9/2.
\end{align*}
If $\lambda_L>0$, dropping the negative density term, using $n^2-1\geq n-1$, and the tangent bound on $e^{2u}$ give
$q_0(L+u)\leq D_Le^{-\lambda_Lu}$ for $u\geq0$.
Consequently for $T_L=\min(|X|,L)$,
\[
 0\leq\E|X|^k-\E T_L^k\leq
 \beta_k(L):=2D_L\sum_{j=1}^k\binom{k}{j}L^{k-j}
                    \frac{j!}{\lambda_L^{j+1}}.
\]
At $L=0.8$, valid rounded bounds are
$\beta_2<3.505\cdot10^{-6}$,
$\beta_4<4.948\cdot10^{-6}$,
$\beta_6<5.289\cdot10^{-6}$.

For $k=2,4,6$, let $\bar V_k,s_k^2$ be the sample mean and unbiased variance of $T_L^k$ and define
\[
 e_k=\sqrt{2s_k^2\log(12/\alpha_m)/N}
       +\frac{7L^k\log(12/\alpha_m)}{3(N-1)}.
\]
With probability at least $1-\alpha_m$, all three target moments lie in
$[\ell_k,u_k]=[\max(0,\bar V_k-e_k),\bar V_k+e_k+\beta_k]$.
Let $d_k=\max(|\ell_k-\theta_{0,k}|,|u_k-\theta_{0,k}|)$. Equation~\eqref{eq:remainder} gives
\[
 R_{\max}=d_2d_6+10|a_0|d_2d_4+5|b_0|d_2^2
             +5d_2^2d_4+\tfrac{10}{3}d_4^2.
\]
Combining this event with the one-sided scalar empirical-Bernstein event proves \eqref{eq:Bbound}. The scalar critical points follow by solving
$g_2+2g_4x^2+3g_6x^4=0$; for the stated anchor their squared positions are approximately $0.08562773$ and $0.33356936$. This verifies monotonicity past $L=0.8$ and the stated range enclosure.

If the actual sampling CDF differs from the target by at most $\varepsilon$, bounded-variation integration by parts bounds the bias of any clipped observable $f$ by $\operatorname{TV}(f)\varepsilon$. For $T_L^k$ the variation is $2L^k$; for the clipped scalar a conservative bound is $6(0.002944)$. Add these deterministic errors to the moment and scalar radii and recompute $R_{\max}$. This transfer does not multiply the law discrepancy by $N$. It controls expectation bias for bounded observables, whereas the concentration theorem still requires independent draws from the actual numerical law.

\end{document}